\documentclass[reqno]{article}
\usepackage{amssymb, amsmath, amsthm, amsfonts, amscd}
\usepackage[english]{babel}

\usepackage{mathtools}
\usepackage{color}
\usepackage{hyperref}
\hypersetup{
    colorlinks=true,
    citecolor=red,
    filecolor=black,
    linkcolor=blue,
    urlcolor=black
}

\chardef\bslash=`\\ 

\newtheorem{thm}{Theorem}[section]
\newtheorem{cor}[thm]{Corollary}
\newtheorem{lem}[thm]{Lemma}
\newtheorem{prop}[thm]{Proposition}
\newtheorem{cond}{Condition}
\newtheorem{conj}{Conjecture}[section]
\newtheorem{notation}{Notation}
\newtheorem{observation}{Observation}

\newtheorem{defn}{Definition}[section]

\newcommand{\N}{\mathbb{N}}

\newcommand{\HH}{\mathcal{H}}

\newcommand{\PP}{\mathcal{P}}

\newcommand{\RR}{\mathcal{R}}

\newcommand{\Z}{\mathbb{Z}}
\newcommand{\Q}{\mathbb{Q}}
\newcommand{\R}{\mathbb{R}}
\newcommand{\C}{\mathbb{C}}
\newcommand{\Sp}{\mathbb{S}}
\newcommand{\T}{\mathbb{T}}
\newcommand{\TT}{\mathcal{T}}

\newcommand{\WW}{\mathcal{W}}
\newcommand{\NN}{\mathcal{N}}
\newcommand{\EE}{\mathcal{E}}

\newcommand{\Diff}{\mathrm{Diff}}

\newcommand{\Id}{\mathrm{Id}}

\def\a{\alpha }
\def\sm{C^{\infty} }
\def\smm{C^{\infty} _{\mu} }

\def\l{\lambda }

\def\n{\nu }
\def\p{\pi}
\def\r{\rho}
\def\s{\sigma}
\def\t{\tau}

\def\d{\delta}

\def\z{\zeta}
\def\w{\omega}
\def\e{\varepsilon}
\def\f{\varphi}
\def\.{\cdot }
\def\ra{\rightarrow}
\def\hra{\hookrightarrow}

\def\begeq{\begin{equation*}}
\def\endeq{\end{equation*}}

\newcommand{\vertiii}[1]{{\left\vert\kern-0.25ex\left\vert\kern-0.25ex\left\vert #1 
    \right\vert\kern-0.25ex\right\vert\kern-0.25ex\right\vert}}

\title{
\textsc{\textbf{Invariant Distributions and local theory of
quasiperiodic cocycles in $\T ^{d} \times SU(2)$}}\\
\author{N. Karaliolios}
\author{Nikolaos Karaliolios \footnote{Universit\'{e} de Lille, LABEX CEMPI
(ANR-11-LABX-0007-01). Email: nikolaos.karaliolios@univ-lille.fr. This work was
supported by the ERC Starting Grant "Quasiperiodic".}}
}

\begin{document}

\nocite{HouPop13}
\nocite{HY2012}

\maketitle

\begin{abstract}
We study the linear cohomological equation in the smooth category
over quasi-periodic cocycles in
$\mathbb{T} ^{d} \times SU(2)$. We prove that, under a full measure condition on the
rotation in $\mathbb{T} ^{d}$, for a generic cocycle in an open set of cocycles,
the equation admits a solution for
a dense set of functions on $\mathbb{T} ^{d} \times SU(2)$ of zero average with respect
to the Haar measure. This property is known as Distributional Unique
Ergodicity (DUE).

We then show that given such a cocycle,  for a generic function
no such solution exists. We thus confirm in this context a conjecture by A. Katok
stating that the only dynamical systems for which the linear cohomological
equation admits a smooth solution for all $0$-average functions with respect
to a smooth volume are Diophantine rotations in tori.

The proof is based on a careful analysis of the
K.A.M. scheme of \cite{KrikAst} and \cite{NKPhD}, inspired by
\cite{El2002a}, which also gives a proof of the local density of
cocycles which are reducible via finitely differentiable or measurable
transfer functions.

\textbf{Key words}: K.A.M. theory, Invariant Distributions, Renormalization, Quasiperiodic cocycles, SU(2)

\textbf{Mathematics Subject Classification}: 37C55
\end{abstract}  

\tableofcontents

\section{Introduction}

\subsection{The setting}

Given a smooth diffeomorphism $f$ on a smooth compact manifold $M$,
the linear cohomological equation over $f$ is
\begin{equation} \label{eq lin cohom}
\psi \circ f - \psi = \f
\end{equation}
where the known function is $\f \in \sm (M)$, and the unknown function is $\psi$.
The operator $\psi \mapsto \psi \circ f - \psi $ is referred to as the \textit{coboundary
operator}, and its image in $\sm$ as the space of \textit{coboundaries over} $f$.
The solution can be sought in $\sm (M)$, in spaces of lower regularity or even
in spaces of distributions (we then speak of weak solutions). Cohomological equations
arise naturally in the
study of Dynamical Systems in various contexts, such as the construction of invariant
tori in K.A.M. theory, the construction of smooth invariant volume forms, or the
study of deviations of ergodic sums from their averages. The linear cohomological
equation, eq. \ref{eq lin cohom}, is the simplest one of the kind, but its study is already
difficult and important.

For example, an important factor of the strength of K.A.M. theory comes from the
fact that for minimal rotations in tori $R_{\a} : \T ^{d} \ni x \mapsto x+\a \in \T ^{d}$, the
image of the corresponding coboundary operator $\psi \mapsto \psi \circ R_{\a} - \psi$ is
of codimension $1$ in $\sm (\T ^{d} )$. This property is called
\textit{Distributional Unique Ergodicity} (DUE) and refines Unique Ergodicity in a non-trivial
manner, since the mapping of $\T ^{2}$ given by
\begeq
(x,y) \mapsto (x+\a , x+y)
\endeq
where $\a \in \R \setminus \Q$ is uniquely ergodic in the classical sense, even though it
is not DUE (proof is by direct calculation, or consult \cite{Kat01}). In the primordial
example of DUE, that of minimal translations in tori
\begeq
R_{\a }: x \mapsto x+ \a \in \T ^{d}
\endeq
we can already distinguish two important subcaces: Diophantine and Liouville translations
(cf. \S \ref{Arithmetics, continued fraction expansion} for the definitions). Diophantine
translations are those given by vectors $\a$ that are approximable only at a polynomial
rate by rational vectors, while Liouville ones can be approximated at an exponential rate
(the rate is with respect to the denominator of the coordinates of the rational vector).
Direct application of the Fourier transformation shows that the range of the coboundary
operator is closed iff the vector $\a$ is Diophantine. If $\a$ is Liouville, then
the range is dense in the space of functions of $0$ average, but a generic function
is not a coboundary. The property of a diffeomorphism $f$ having a closed space of
coboundaries is known as \textit{Cohomological Stability}.

Only recently, cf. \cite{AFKo2015}, were
constructed examples of manifolds other than tori supporting DUE diffeomorphisms, thus providing
counterexamples to a conjecture by G. Forni (see
\cite{FlFoRH13}). This conjecture concerned the spaces that carry
DUE $C^{\infty}$-diffeomorphisms, and stated that such spaces,
if they are closed manifolds, should be tori. One of the classes of counterexamples constructed was a subset of
$SW ^{\infty} (\T , P) $, the space of skew-product diffeomorphisms of $\T \times P$,
with $P $ a homogeneous space of the compact type, i.e. a manifold $P=G/H$ where
$G $ is a compact Lie group and $H$ is a subgroup. Skew-product diffeomorphisms are
dynamical systems on $\T ^{d} \times P$ of the type
\begeq
(\a , A(\. )) : (x,p) \mapsto  (x+ \a , A(x) .p)
\endeq
where $p = s.H \in P$, and $s \in G$. The mapping $A(\. ): \T ^{d} \ra G$ is $\sm $ smooth.
The space of such dynamical systems, called \textit{quasi-periodic cocycles over} $\a$
(or simply cocycles) is denoted by $SW^{\infty }_{\a} (\T ^{d} , P)$ when
$\a \in \T ^{d}$, called the frequency, is considered fixed. We call
\begeq
SW^{\infty } (\T ^{d} , G) = \cup _{\a \in \T ^{d}} SW^{\infty }_{\a} (\T ^{d} , G)
\endeq
Conjugation in this space is fibered conjugation. If we let $B(\. ): \T ^{d} \ra G$ be a
mapping whose regularity can range from measurable to $\sm$-smooth, we can change coordinates
fiberwise in $\T ^{d} \times P $ by
\begeq
(x,p) \mapsto (x, B(x).p)
\endeq
which results in conjugating the cocycle $(\a , A(\. ))$ to
\begeq
(\a , \tilde{A}(\. )) = ( \a ,B(\. + \a). A(\. ).B^{-1}(\. ))
\endeq
(cf. \S \ref{sec def dynamics} for more details).

In \cite{AFKo2015}, the authors proved that DUE is generic in the closure of
\begeq
\{ ( \a ,B(\. + \a).\Id . B^{-1}(\. )) , \a \in \Q , B(\. ) \in \sm (\T , G) \} \subset
SW^{\infty} (\T , P)
\endeq
thus letting the frequency $\a$ vary. It is a classical fact that Liouville numbers,
those whose coboundary space is not closed, form a residual set in $\T$, and as a
consequence the authors of \cite{AFKo2015} were not able to address cohomological
stability, due to its obvious failure for observables depending only on the variable in
$\T $.

Our goal is to strengthen their result, more specifically for $G=SU(2)$, on the basis of the recent development
of the theory of quasi-periodic cocycles in compact Lie groups in our PhD thesis \cite{NKPhD}.
We thus obtain a theorem of a different flavour than that of \cite{AFKo2015},
where cocycles over a fixed Diophantine frequency are considered. As a result, the obvious
obstruction to Cohomological Stability mentioned here above is lifted. On the other hand,
the algebraic setting is more
restricted in order to keep the arguments transparent, even though we are confident that
our results can be generalized to all semisimple compact Lie groups. We thus prove
genericity of DUE in an open subset of $SW^{\infty} _{\a} (\T , SU(2))$, with $\a $
satisfying a Diophantine condition.

This open subset is called the K.A.M. regime for such cocycles, for the reason that the
Diophantine condition on the frequency, coupled with the assumption that the mapping
$A( \. )$ defining the cocycle is close to a constant mapping, allow the application of
the K.A.M. machinery. The resulting theorem is known as the Almost Reducibility theorem
(cf. \cite{El2001}, \cite{KrikAst}, \cite{NKPhD}), since it establishes that every cocycle
in this open set can be conjugated arbitrarily close to constant ones.

This is carried out by a K.A.M. scheme whose construction is as follows (cf. \S \ref{sec a.r. and due}
for the details, or \cite{NKContSpec} for a description). Assume that
$\a $ satisfies a Diophantine condition and that the mapping $A(\. ) = A.e^{F(\. )}$ is
close to the constant $A \in SU(2)$ (which amounts to $F(\. ): \T ^{d} \ra su (2)$, the
Lie algebra of $SU(2)$, being close to the $0$ mapping). Then, an iterative scheme can be
constructed, conjugating sucessively the cocycle
\begeq
(\a , A.e^{F(\. )}) = (\a , A_{1}.e^{F_{1}(\. )})
\endeq
to the cocycle $(\a , A_{n}.e^{F_{n}(\. )})$, where $F_{n}(\. ) \overset{\sm }{\ra } 0 $.
This is precisely the definition of almost reducibility. The conjugation is constructed in
two steps, at each iteration of the scheme. At the first step, a cohomological equation
is solved (a \textit{twisted} one, slightly more complicated than the linear one) for a finite
number of Fourier modes of the perturbation. If the eigenvalues of $A_{n}$ are away
from $\{ e^{2i \pi k \. \a} \} _{k \in \Z ^{d}}$ (we then say that $A _{n}$ is \textit{non-resonant}), then the
equation can be solved with good estimates. The solution provides a close-to-the-$\Id$
conjugation reducing the perturbation $F_{n}$ to a second order one, and the iteration
is complete. In this case, the commutator between $A_{n}$ and $A_{n+1}$ is of the order
of $F_{n}$.
If $A_{n}$ is resonant, i.e. if its eigenvalues are very close to some $e^{2i \pi k _{n} \. \a}$, then the
perturbation is still of the first order, but with only one Fourier mode active, the mode
$k_{n}$, called the resonant mode. The perturbation is then reduced to a second-order one
via a far-from-the-$\Id$ conjugation, and the iteration is then complete, but in this case
the commutator between $A_{n}$ and $A_{n+1}$ can be of any order.

The fact that this theorem comes with very precise estimates on the norms of the conjugations and the
size of the $F_{n} (\. )$ allows for a very efficient analysis of the dynamics for the
cocycles in the K.A.M. regime. The study of the cohomological equation over a generic
cocycle is actually made possible thanks to this analysis, but let us first mention some
more direct implications.

A cocycle is called reducible if it can be conjugated to a constant one, i.e. a cocycle of
the form $(\a , A)$ with $A \in G$, and the regularity of the conjugation is to be made
precise in each statement. A direct corollary of the Almost Reducibility theorem is the
density of $\sm$ reducible cocycles in the K.A.M. regime. Another consequence of almost
reducibility, more difficult to obtain, is the ergodicity of a generic cocycle in the K.A.M.
regime for the Haar measure on $\T ^{d} \times G$. This was obtained by H. Eliasson
(\cite{El2001}), by imposing a condition on the commutator between $A_{n}$ and $A_{n+1}$.
The condition, put informally, imposes that the commutator stays bounded away from
the $\Id$, along a subsequence of K.A.M. iterations.

In the present article we refine and strengthen his results. Namely, we show that a cocycle
is reducible with a conjugation in a given Sobolev space $H^{\s}$ if, and only if, the sequence of
such commutators (or more precisely their distance to the $\Id$)
is summable in the corresponding Sobolev space $h^{\s }$.
We also show that a cocycle will be DUE if, and only if, this commutator takes values
infinitely often sufficiently close to a countable set.

In a subsequent article, \cite{NKContSpec}, we pushed the techniques of the present paper
further. We applied them in the solution of the eigenfunction equation,
\begeq
\f \circ (\a , A(\. )) = \l \f
\endeq
instead of the linear cohomological one, and obtained that for a generic cocycle in the
K.A.M. regime, no $L^{2}$ solution of this equation exists. This implies that the
cocycle is weak mixing in the space orthogonal to functions depending only in the variable
in $\T ^{d}$. We encourage the reader to consult that article, because the setting
therein is more classical.

More recently in \cite{NKCohomRigAK}, we were able to establish that, for general reasons,
the DUE cocycles obtained in the present article are not Cohomologically stable due to
the fact that they are exponentially fast approximated by resonant cocycles.
We were able to prove that no DUE and Cohomologically stable cocycles
exist in the K.A.M. regime of $SW^{\infty}_{\a} (\T ^{d} ,P)$ with $P$ any
homogeneous space of compact type (instead of $SU(2)$ as herein). The
description of the dynamics in $SW^{\infty}_{\a} (\T ^{d} ,P)$ remains,
however, far less precise than that of $SW^{\infty}_{\a} (\T ^{d} ,SU(2))$.


\subsection{Statement of the results}

We will begin by recalling what we need from the theory developed in \cite{NKPhD},
inspired by \cite{KrikAst} and
\cite{El2002a}, translated in the context of the article. It is the
local theory (\S 8 of \cite{NKPhD}), and more precisely the proof of the theorem
on local almost reducibility. Then, we show that if we let the
condition devised by H. Eliasson in \cite{El2002a} degenerate in a
controlled way, we obtain cocycles that are reducible in any given
Sobolev regularity. The proof of this theorem occupies \S \ref{loc dens of fin dif}.
\begin{thm} \label{finite reduc thm}
Let $\a \in DC(\gamma , \tau ) \subset \T ^{d}$, and let
$\WW \subset SW^{\infty}_{\a} (\T ^{d} ,SU(2))$ be the neighborhood
of constant cocycles defined in the local almost reducibility
theorem (see, e.g. \cite{KrikAst} or thm \ref{thm local a.r.} ). Then,
for any given $\s_{0}\in [0,\infty ]$,
$H^{\s_{0}}$-reducible cocycles are $C^{\infty}$-dense in $\WW $.
\end{thm}

By "$H^{\s_{0}}$-reducible cocycle", we mean a $C^{\infty}$-smooth cocycle which
is conjugate to a constant one, and the transfer function is of
Sobolev regularity $H^{\s _{0}}$, but not $H^{s}$, for any $s >\s _{0}$.
We remark that the proof implies the same conclusion for finite
and negative $\s _{0}$, but we ignore its dynamical content. Finally,
since the result is already known for $\s  _{0}= \infty$, see
\cite{KrikAst} or \cite{NKPhD}, we will focus on the remaining case of finite
$\s _{0}$. The whole discussion in paragraphs \ref{remarks and cons} and
\ref{loc dens of fin dif} shows the optimality of the K.A.M. scheme
for cocycles over a diophantine rotation, satisfying a smallness
condition.

The next theorem and the flowchart proving it show that the K.A.M. scheme provides all
the ingredients necessary for the construction of an $H^{\s }$-smooth transfer function
reducing a given cocycle, provided that such a function exists.
\begin{thm}
For cocycles in $\WW$, the property of $H^{\s}$-reducibility, $ \s \in [0, \infty]$, is
algorithmically decidable.
\end{thm}
The algorithm ends after an at most countable number of steps.
\begin{proof}
Let a cocycle in $\WW $.
Then, using the K.A.M. scheme, one can determine whether the cocycle is reducible and
construct a sequence of $C^{\infty}$-smooth transfer functions, converging to the
reducing transfer function in the $H^{\s}$ topology using the following flowchart. We
remind that the series $\theta $ measures the failure of commutativity between the
constants $A_{n_{i}}$ and $A_{n_{i+1}}$ as constructed by the K.A.M. scheme.
\begin{enumerate}
\item If the product of conjugations constructed by the K.A.M. scheme converges, it
converges in $C^{\infty }$ and the cocycle is $C^{\infty }$-reducible (cf. lemma
\ref{lemconv conj smooth}).
\item If not, check the summability in $h^{s}$ of the series $\theta$
defined in eq. \ref{cvgence cdtion}. Then,
\begin{itemize}
\item If there exists $s \geq 0 $ such that $\theta \in h^{s}$ and
$\s \in [0,\infty ] $ is the maximal such index, the cocycle is
$H^{\s }$-reducible (and not in any higher regularity).
\begin{itemize}
\item If $\s = \infty$, calculate the quantity $\bar{\n } $
of the proof of prop. \ref{thm cohom stab implies reduc} (eq.
\ref{eq def nbar}). If
$\bar{\n } < \infty $, change the parameters of the K.A.M. scheme
as in the proof and define a scheme for which the product of
conjugations converges. This is always the case for cocycles measurably
conjugate to a Diophantine constant, by \cite{NKRigidity}.
\item If $\s = \infty$ and $\bar{\n } = \infty $, or if $\s < \infty$,
run the algorithm of lemma \ref{conj in fin reg iter} and construct a
reducing conjugation. The constant to which the cocycle is reduced is
Liouville, by \cite{NKRigidity}.
\end{itemize}
\item If no such $s $ exists, the cocycle is not $L^{2}$ reducible, since it is
weak mixing in the fibers by \cite{NKContSpec}.
\end{itemize}
\end{enumerate}
This concludes the proof.
\end{proof}


Subsequently, we use the K.A.M. scheme in order to solve the
cohomological equation in $C^{\infty } (\T ^{d} , SU(2))$, $d\geq 1$,
up to an arbitrarily small error. This is equivalent to Distributional
Unique Ergodicity of the cocycle. It can be done under a generic
condition on the cocycle, inspired by the one in \cite{El2002a}. Under this condition,
for any given function of $0$ average, if we go deep
enough in the K.A.M. scheme that reduces the given cocycle and solve the cohomological
equation over the cocycle $(\a , A_{n})$ in the corresponding system of coordinates,
the error that we commit becomes arbitrarily small in
the $C^{\infty}$ topology.
\begin{thm} \label{thm DUE}
Let $\a \in DC(\gamma ,\tau) \subset \T ^{d} $ with $d \in \N ^{*} $.
Then, DUE holds in a $G_{\delta}$ dense subset of $\WW$.
\end{thm}
This is a local genericity theorem, and the generic condition that
is needed for DUE is strictly stricter
than the one of \cite{El2002a}. We actually show the optimality of the condition, in the sense that if
we let it degenerate so that the cocyle thus obtained is Uniquely Ergodic,
functions of $0$ average exist for which the equation cannot be
solved approximately (see prop. \ref{optim of cdtion}).


These theorems hold for cocycles of arbitrary numbers of frequencies.
The following corollary uses the global density theorem
\ref{thm global density}, which has been proved only for one-frequency
cocycles (\cite{Krik2001},\cite{Fra2004},\cite{NKPhD}).
\begin{cor} \label{cor DUE glob}
If $\a \in RDC \subset \T $, DUE is generic in $SW^{\infty}_{\a } (\T ,SU(2))$.
\end{cor}
\begin{proof}
The density part of the corollary is obtained in the same way as
corollary \ref{global density of UE} from theorem \ref{thm global density}.
Moreover, DUE is a $G_{\d }$ condition in the $C^{\infty}$ topology
by lemma \ref{lem DUE Gd}.
\end{proof}
This was proved in \cite{AFKo2015} for cocycles over generic
Liouvillean rotations, but it is also found under this very
restrictive (in the topological sense) arithmetic condition. In the same
way we can obtain
\begin{cor}
DUE is generic in $SW^{\infty}(\T , SU(2))$.
\end{cor}
This is, of course, a special case of theorem \ref{DUE AFKo}. The argument goes like the one in the proof of
corollary \ref{genericity in hom to csts}.
\begin{proof}
DUE is a $G_{\d }$ property by lemma \ref{lem DUE Gd}. By corollary \ref{cor DUE glob},
it holds in a dense subset of $SW_{\a }^{\infty}(\T , SU(2))$, for a dense subset of $\a \in \T$.
\end{proof}

In the results above $SW^{\infty}_{\a } (\T ,SU(2))$ can be replaced by $SW^{\infty ,1}_{\a } (\T ,SO(3))$,
the space of cocycles in $\T \times SO(3) $ homotopic to the $\Id $. On the other hand, DUE seems to be a rare
property in the remaining connected component of $SW^{\infty }_{\a } (\T ,SO(3))$.

Finally, we address the problem of cohomological stability and prove the
following rigidity theorem.
\begin{prop} \label{thm cohom stab implies reduc}
If a cocycle in $\WW$ is cohomologically stable, it is reducible.
\end{prop}
In fact, cohomological stability implies a very good control of the conjugations produced by the K.A.M. scheme:
the latter can be made to produce a conjugation by adjusting the
parameter measuring the distance from resonant constants (the parameter $\n$),
so that the cocycle
is actually reducible (to a Diophantine or resonant constant).
Since the remaining part of the total space $SW_{\a}^{\infty}(\T , SU(2)), \a \in RDC$,
is filled by cocycles conjugate to diagonal ones,
we obtain the following result.
\begin{cor} \label{no cout-ex}
If $\a \in RDC$, there are no counterexamples to conjecture
\ref{conj Katok} in $SW^{\infty}_{\a }(\T , SU(2))$.
\end{cor}

\begin{proof}
By proposition \ref{thm cohom stab implies reduc}, an almost reducible cohomologically
stable cocycle is reducible. Reducible cocycles preserve a foliation of
$\T \times SU(2)$ into tori $\T ^{1} \times \Sp ^{1}$ and are thus not DUE.

By theorem \ref{normal form SU(2)}, cocycles that are not almost reducible preserve a
foliation into tori $\T ^{1} \times \Sp ^{1}$, and are thus not DUE.
\end{proof}

%

A detailed classification of the dynamics in $\WW $ based on theorems proved herein, in the
bibliography, as well as in the reference itself, is given in \cite{NKContSpec}, \S 2.

\textbf{Acknowledgments} We would like to thank A. Avila, B.
Fayad and A. Kocsard for drawing our attention to the problem of DUE,
and H. Eliasson for his interest in the author's work during the last
years and for asking the question that lead to thm.
\ref{finite reduc thm}. The author would also like to thank
Rapha\"{e}l Krikorian for his support during the preparation of the
author's PhD thesis and his subsequent works, as well as the referee for their
remarks.

\section{Facts from algebra and arithmetics}

\subsection{The group $SU(2)$}
The matrix group $SU(2)$ is the multiplicative group of unitary
$2 \times 2$ matrices of determinant $1$. For the greatest
part of the article, unless the contrary is explicitly stated,
$G$ will stand for the group $ SU(2) $.

Let us denote the matrix $S\in G$, $S=\begin{pmatrix}
z & w \\ 
-\bar{z} & \bar{w}
\end{pmatrix}
$, where $(z,w)\in \C^{2}$ and $|z|^{2}+|w|^{2}=1$, by $\{z,w\}_{G }$. The subscript will be
suppressed from the notation, unless necessary. The manifold $G=SU(2)$ is thus naturally identified with
$\Sp ^{3}\subset \C^{2}$, and, in particular, $G$ is simply connected. The law of multiplication is pushed
forward to the mapping of $\Sp ^{3}\times \Sp ^{3}\ra \Sp ^{3}$ given by
$\{z_{1},w_{1}\}.\{z_{2},w_{2}\}=\{z_{1}z_{2}-w_{1}\bar{w}_{2},z_{1}w_{2}+\bar{z}_{2}w_{1}\}$.  Inversion is pushed
forward to the involution of $\Sp ^{3}$, $\{z,w\}\mapsto \{\bar{z},-w\}$. When coordinates in $\C ^{2}$ are fixed, the
circle $\TT \equiv \Sp ^{1}$ is naturally embedded in $G$ as the group of diagonal matrices, which is a maximal torus (i.e. a
maximal abelian subgroup) of $G$.

The Lie algebra $g=su(2)$ is naturally isomorphic to $\R^3 \approx \R \times \C$
equipped with its vector and scalar product.
It will be denoted by $g$. The element $s=\begin{bmatrix}
it & u \\ 
-\bar{u} & -it
\end{bmatrix}
$ will be denoted by $ \{t,u\}_{g} \in \R \times \C$. The scalar product is defined by
\begin{equation*}
\langle \{t_{1} ,u_{1} \} ,\{t_{2} ,u_{2} \} \rangle =
t_{1}t_{2}+\mathcal{R} (u_{1}\bar{u}_{2})=t_{1}t_{2}+\mathcal{R} u_{1}.\mathcal{R} u_{2}+
\mathcal{I} u_{1}.\mathcal{I} u_{2}
\end{equation*}
Mappings with values in $su(2)$ will be denoted by
\begin{equation*}
U(\. ) = \{ U_{t}(\. ), U_{z}(\. ) \}_{g}
\end{equation*}
in these coordinates, where $U_{t}(\. )$ is a real-valued and $U_{z}(\. )$ is a complex-valued function.

The adjoint action of $h\in su(2)$ on $su(2)$ is pushed-forward to twice the vector product:
\begin{equation*}
ad_{ \{ 1,0 \} }.\{ 0,1 \} =  [ \{1,0\}, \{0,1\} ]  = 2\{ 0,i \}
\end{equation*}
plus cyclic permutations, and the Cartan-Killing form, normalized by
$\left\langle h,h^{\prime }\right\rangle =-\frac{1}{8\pi }tr(ad(h)\circ ad(h^{\prime }))$ is pushed-forward to the scalar
product of $\R^{3}$. The periodic geodesics of the group for the induced
Riemannian structure are of the form $S.E_{r}(\. ).S^{-1} =  S.\exp (\{2\pi r\.  ,0\}_{g}).S^{-1}$, $S\in G$ and $r \in \N ^{*}$.

Under this normalization, the minimal length for a geodesic
leaving from the $\Id$ and arriving at $-\Id$ is $1/2 $, and the minimal
length of a periodic geodesic is twice as much. We also find
directly that the preimages of the $\Id$ in the maximal toral (i.e. abelian) algebra of
diagonal matrices are points of coordinates in the lattice $2\pi \Z$.

The adjoint action of the group on its algebra is pushed-forward to the
action of $SO(3)$ on $\R \times \C$. In particular, the diagonal matrices, of the form $S = \exp (\{2\pi s,0\}_{g})$,
$Ad(S).\{ t,u\} =\{ t,e^{4i\pi s }u\} $.

Finally, the Weyl group of $SU(2)$ (i.e. the group of symmetries of
$\TT \hra SU(2) $ induced by $SU(2)$) is naturally isomorphic to the
multiplicative group $\{1,-1\}$ and two representatives are the matrices%
\begin{equation*}
\Id \text{ and }
\begin{bmatrix}
0 & 1 \\ 
-1 & 0
\end{bmatrix}
\end{equation*}
the action of the second one being a reflection with respect to the complex
plane:
\begin{equation*}
(-1)_{W}.\{t,0\}_{su(2)}=\{-t,0\}_{su(2)}
\end{equation*}
The \textit{normalizer} $ \NN $ of the subgroup $\TT $ of diagonal matrices contains $\TT $ itself, and matrices of the form
$\begin{pmatrix}
0 & 1 \\
-1 & 0
\end{pmatrix}
. \TT $. The diagonalization of matrices in $G$ uses elements of $G \mod \NN $.

\subsection{Functional Spaces}

We will consider the space $C^{\infty }(\T ,g)$ equipped with the
standard maximum and the Sobolev norms. For $s \geq 0$, we define
\begin{eqnarray*}
\left\Vert U\right\Vert _{s} &= & \max_{0\leq \sigma \leq s} \max_{\T }\left\vert \partial
^{\sigma }U(\. )\right\vert \\
\left\Vert U\right\Vert _{H^{s}}^{2} &= & \sum_{k\in \Z^{d}}(1+|k|^{2})^{s}|\hat{U}(k)|^{2}
\end{eqnarray*}
where $\hat{U}(k)=\int U(\. )e^{-2i\pi kx}$ are the Fourier coefficients
of $U(\. )$. The injections $H^{s +1}(\T ,g) \hra C^{s}(\T ,g) $
and $C^{s}(\T ,g) \hra H^{s}(\T ,g)$ for all $s \geq 0$ are
continuous .

By abusing the notation, we note $H^{0} = L^{2}$, and define the spaces
\begin{equation*}
\tilde{H}^{\s} = H^{\s} \setminus \{ \cup _{s >\s} H^{s } \}, \s \geq 0
\end{equation*}
The spaces with tildes are generic subspaces of the
corresponding Hilbert spaces, and consist of the elements that
exhibit no additional Sobolev regularity. We will denote the
corresponding spaces of complex sequences by lowercase letters,
\begin{equation}\label{def h spaces}
\begin{array}{r@{}l}
h ^{s} &= \{ f\in \ell ^{2}, \sum (1+n)^{2s}|f_{n}| ^{2} <
 \infty\} \\
\tilde{h} ^{\s} &= h ^{\s} \setminus \{ \cup _{s >\s}
h ^{s } \}
\end{array}
\end{equation}

For this part, see \cite{FollandHarmAn} and \cite{SteinWeissFourierEuc}. In view of the identification of the $C^{\infty}$ manifold $G = SU(2)$ with
$\Sp ^{3} \subset \C ^{2}$, equipped with the natural measure normalized so that sphere has volume $1$,
the space $C^{\infty } ( G  )$ of smooth $\C$-valued functions defined on $ G $, can be
identified with $C^{\infty } ( \Sp ^{3} )$, and the identification is an isometry between the $L^{2}$ spaces.
In order to obtain a basis facilitating the calculations, we use the
representation theory of $G$. Given a system of coordinates $(\z ,\w )$
in $\C ^{2}$, we can define an orthonormal basis for $\PP _{m}$,
the space of homogenous polynomials of degree $m$, by $\{ \psi _{l,m} \}_{0\leq l \leq m}$ where
$ \psi _{l,m} (\z ,\w ) = \sqrt{\frac{(m+1)!}{l!(m-l)!}}  \z ^{l} \w ^{m-l}$. The group $G$ acts on $\PP _{m}$ by
\begeq
\{ z,w \}.\phi (\z ,\w ) = \phi (z \z + w \w , -\bar{w} \z + \bar{z} \w )
\endeq
and the resulting representation is noted by $\p _{m}$.
For $m$ fixed, we can define the matrix coefficients relative to the basis by
\begin{equation} \label{eq def pi}
\p _{m}^{j,p} \{ z,\bar{z} , w , \bar{w}  \} \mapsto \langle \{ z,w \}.\psi _{j,m} ,\psi _{p,m} \rangle
\end{equation}
The matrix coefficients are
harmonic functions of $z,\bar{z},w,\bar{w}$, and are of bidegree $(m-p,p)$, i.e. they are homogeneous of degree $m-p$ in
$(z,w)$, and homogeneous of degree $p$ in $(\bar{z},\bar{w})$, and they generate the space $\EE _{\p _{m}}$.

We thus obtain the decomposition
$L^{2} = \oplus _{m \in \N} \EE _{\p _{m}} = \oplus _{m \in \N } \oplus _{0 \leq p \leq m}  \HH _{m,p}$, where $\HH _{m,p}$
is spanned by $\{ \pi ^{j,p}_{m} , 0\leq j \leq m \}$. In this
decomposition, the harmonic functions are regrouped according to
their bidegree.

Therefore, given a system of coordinates in $\C ^{2}$, a function $f \in L^{2}( \Sp ^{3} ) $
can be written in the form
\begin{equation*}
f(z,\bar{z},w , \bar{w}) = \sum _{m \in \N } \sum_{0 \leq p \leq m} \sum _{0 \leq j \leq m } f^{m}_{j,p }
\p _{m}^{j,p} (z,\bar{z},w , \bar{w})
\end{equation*}
where $f^{m}_{j,p } \in \C $ are the Fourier coefficients.
The functions $\p _{m}^{j,p} (z,\bar{z},w , \bar{w})$ are the eigenvectors of the Laplacian
on $\Sp ^{3}$ and consequently smooth (in fact real analytic), and they form an orthonormal basis for
$L^{2} (\Sp ^{3})$. In higher regularity, they generate a dense subspace of $C^{\infty}$.

The group $G $ acts on $C^{\infty } ( G  ) \equiv C^{\infty } (\Sp ^{3} ) $ by pullback: if $A \in G $ and
$\begin{pmatrix}
z \\
w
\end{pmatrix} \in \Sp ^{3} $, then, for $\phi : \Sp ^{3} \ra \C $,
\begin{equation*}
(A.\phi) 
\begin{pmatrix}
z \\
w
\end{pmatrix}
= \phi
\left( A^{*}
\begin{pmatrix}
z \\
w
\end{pmatrix}
\right)
\end{equation*}
If coordinates are chosen so that $A = \{ e^{2i\pi a } , 0 \} $ is diagonal, then
\begin{equation*}
A.\phi (z,w) = \phi( e^{-2i\pi a }z , e^{2i\pi a }w)
\end{equation*}
and $A$ then acts on harmonics by
\begin{equation*}
A.\p _{m}^{j,p} (z,\bar{z},w , \bar{w}) = e^{2i\pi (m-2p) a } \p _{m}^{j,p} (z,\bar{z},w , \bar{w})
\end{equation*}
where, of course, $m -2p = m-p-p$ is the difference of the degrees of homogeneity in $(\bar{z}, \bar{w})$
and $(z,w)$. Therefore, the harmonics in these coordinates are eigenvectors for the associated operator.
In particular, if $a $ is irrational, the eigenvectors for the eigenvalue $1$ are exactly the elements
$\p ^{j,m/2}_{m}$, $0 \leq j \leq m$.

The group of symmetries
of $\C . \psi _{m/2,m} $ is exactly $\NN $, the normalizer of $\TT $. The subspaces $\C . \psi _{l,m} , l \neq m/2 $
appear to have a smaller normalizer, $\TT $, but the normalizer of $\C \psi _{l,m}  + \C \psi _{m-l,m} $ is
indeed $\NN $. The following two lemmata concern the effect of a change of basis in $\C ^{2}$ on the
Fourier development of $f$ and are proofs and quantifications of this fact.
\begin{lem} \label{change of coord harm}
Let $\psi _{l,m}  \in \PP _{m} $. Let also $\pi _{p,m} $ the projection on $\C \, \psi _{p,m}  ( \z , \w )$.
Then, unless $l' \in \{ l,m-l \}$, there exists a constant $C_{l,m} < 1$ such that for every $A \in G$
\begin{equation*}
\| \pi _{l' ,m} (A. \psi _{l,m} ( \z , \w ) )) \|  \leq C_{l,m}
\end{equation*}
\end{lem}

In fact we show that the only possible way of transforming a monomial in $\PP _{m}$ to another one by a linear
change of coordinates is by permuting the coordinates (and changing the argument within each complex plane),
so that $\C \psi_{l,m} ( \z , \w ) + \C \psi_{m-l,m} ( \z , \w )$ is preserved exactly by $\NN $.

\begin{proof}
Elementary calculation shows that
\begin{equation*}
| \psi _{l,m} ( \z , \w ) | \leq  \sqrt{\frac{(m+1)!}{l!(m-l)!}} \left(  \frac{l}{m} \right) ^{l/2}
\left( \frac{m-l}{m} \right) ^{(m-l)/2}
\end{equation*}
and, if $l<m$, the locus where the maximum is attained is
\begin{equation*}
\{ ( \sqrt{\frac{l}{m}}e^{2i\pi \theta} , \sqrt{\frac{m-l}{m}} e^{-2i\pi \frac{m}{m-l} \theta } ) \}_{\theta \in \T }
\end{equation*}
If $l=m$, the maximum is is attained at $\{ (e^{2i\pi \theta} , 0) \}_{\theta \in \T }$. Clearly, the loci
corresponding to harmonics not obtained by permutations of coordinates cannot obtained by linear transformations of
the sphere, and therefore the only vector in the standard basis of $\PP _{m}$ that can be obtained by the action
of $G$ on $\psi _{l,m}$ is $\psi _{m-l,m} ( \z , \w )$.
In particular, the only transformations leaving $\C \psi _{m/2,m} ( \z , \w )$ fixed are
\begin{equation*}
\begin{pmatrix}
e^{2i\pi \theta} & 0 \\
0 & e^{-2i\pi \theta}
\end{pmatrix}
\text{ and }
\begin{pmatrix}
0 & e^{2i\pi \theta} \\
- e^{-2i\pi \theta} & 0
\end{pmatrix}
\end{equation*}
for $\theta \in \T $.
\end{proof}

The second lemma examines the effect of changes of coordinates on the eigenvectors for the eigenvalue $1$,
$\psi _{m/2,m} $.
\begin{lem} \label{Legendre roots}
For a given $m >0 $ and even, there exist $m$ points $D \in G / \NN $ such that
$\pi _{m/2,m} (D.\psi _{m/2,m}) = 0$.
\end{lem}
\begin{proof}
Call $l=m/2$. We can calculate explicitly the projection:
\begin{equation*}
\pi _{l,m} (\{ z,w \} _{G}.\psi _{l,m}) = \sum _{0}^{l} (-1)^{i}
\begin{pmatrix}
l \\
i
\end{pmatrix} ^{2}
|z|^{2(l-i)}|w|^{2i} \psi _{l,m} = p_{l} (|z|,|w|) \psi _{l,m}
\end{equation*}
The factor of the projection, $p_{l}$, is a Legendre polynomial of the variable $|z|^{2}$
and $|w|^{2}=1-|z|^{2}$. The fact that all its roots are real and in the interval $[0,1]$ is classical,
see for example \cite{ArnPDE}.
\end{proof}
Returning to more general facts from calculus, the $C^{s }$ norms for functions in $C^{\infty } (G )$ are defined in a
classical way, and the Sobolev norms are defined by
\begin{equation*}
\left\Vert f \right\Vert _{H^{s}}^{2}=\sum_{m,j,p }(1+m^{2})^{s}|f_{m}^{j,p}|^{2}
\end{equation*}
where $f_{m}^{j,p} $ are the coefficients of the harmonics in the expansion of $f $.

We will also use the convexity or Hadamard-Kolmogorov inequalities (see \cite{KolmIneq49}) ($U\in
C^{\infty }(\T ,g)$):
\begin{equation}
\left\Vert U(\. )\right\Vert _{\sigma }\leq C _{s, \sigma }\left\Vert
U\right\Vert _{0}^{1-\sigma /s}\left\Vert U\right\Vert _{s}^{\sigma /s}
\label{hadamard}
\end{equation}
for $0\leq \sigma \leq s$, and the inequalities concerning the composition of functions (see \cite{KrikAst}):
\begin{equation} \label{composition}`
\left\Vert \phi \circ (f+u)-\phi \circ f\right\Vert _{s}\leq
C_{s}\left\Vert \phi \right\Vert _{s+1} \vertiii{f}_{s} \left\Vert u\right\Vert _{s}
\end{equation}
where we introduce the notation
\begin{equation*}
\vertiii{f}_{s} =
(1+\left\Vert f\right\Vert
_{0})^{s}(1+\left\Vert f\right\Vert _{s})
\end{equation*}

We will use the truncation operators for mappings $\T \ra g$ defined by
\begin{eqnarray*}
T_{N}f(\. ) &=&\sum_{|k|\leq N}\hat{f}(k)e^{2i\pi k\. } \\
\dot{T}_{N}f(\. ) &=&T_{N}f(\. )-\hat{f}(0) \\
R_{N}f(\. ) &=&\sum_{|k|>N}\hat{f}(k)e^{2i\pi k\. }\text{ and} \\
\dot{T}f(\. ) &=&Tf(\. )-\hat{f}(0)
\end{eqnarray*}
These operators satisfy the estimates
\begin{eqnarray} \label{truncation est}
\left\Vert T_{N}f(\. )\right\Vert _{C^{s}} &\leq
&C_{s}N \left\Vert f(\. )\right\Vert _{C^{s}} \\
\left\Vert R_{N}f(\. )\right\Vert _{C^{s}} &\leq &C_{s,s'} N^{s-s^{\prime }+2} \left\Vert
f(\. )\right\Vert _{C^{s^{\prime }}}
\end{eqnarray}
The Fourier spectrum of a function will be denoted by $\hat{\sigma}(f)=\{k\in \Z,~\hat{f}(k)\not=0\}$.

Finally, for functions $f: G \ra \C $, we truncate on $m$:
\begin{equation*}
T_{M}f(\. ) = \sum_{m \leq M} \sum _{0 \leq j,p \leq m } f^{m}_{j,p} \p _{m}^{j,p}
\end{equation*}
and likewise for the other operators. These operators satisfy the estimates
\begin{eqnarray*}
\left\Vert T_{M}f(\. )\right\Vert _{C^{s}} &\leq
&C_{s}M^{2 }\left\Vert f(\. )\right\Vert _{C^{s}} \\
\left\Vert R_{M}f(\. )\right\Vert _{C^{s}} &\leq &C_{s,s'} M^{s-s^{\prime }+2}\left\Vert
f(\. )\right\Vert _{C^{s^{\prime }}}
\end{eqnarray*}

\subsection{Arithmetics, continued fraction expansion} \label{Arithmetics, continued fraction expansion}

For this section, see \cite{KhinContFr} or \cite{ArnGeomMeth}. Let us
introduce some notation. For $\a \in \R ^{*}$, define  $ |\a | _{\Z } =dist(\a  ,\Z) = \min _{\Z } |\a - l| $
to be the distance from the closest integer. We also denote by $[\a  ]$, resp. $\{\a\}$ the integer part,
resp. the fractional part, of $\a$. We finally call $\mathrm{G}(\a)=\{\a^{-1}\}$, the Gauss map.

Consider $\a\in \T \setminus \Q$ fixed, and let $p_{-2}=q_{-1}=0$ and $p_{-1}=q_{-2}=1$. Then
$(p_{n}/q_{n})_{n\geq 0}$ is the sequence of best rational approximations of $\a$ defined recursively as follows.
Let $\a_{n}=\mathrm{G}^{n}(\a  )=\mathrm{G}(\a_{n-1})$, $a_{n}=[\a_{n-1}^{-1}]$, and
$\beta _{n}=\prod\limits_{0}^{n}\a_{k}$. Then
the Euclidean division of $\beta _{n-2}$ by $\beta _{n-1}$ reads
\begin{equation*}
\beta _{n-2}=a_{n}\beta _{n-1}+\beta _{n}
\end{equation*}
and we can define
\begin{eqnarray*}
q_{n} &=&a_{n}q_{n-1}+q_{n-2} \\
p_{n} &=&a_{n}p_{n-1}+p_{n-2}
\end{eqnarray*}
$p_{n}$ and $q_{n}$ are strictly positive for $n\geq 1$ and $\beta
_{n}=(-1)^{n}(q_{n}\a-p_{n})$. We have
\begin{eqnarray*}
\frac{1}{q_{n}+q_{n+1}} &<&\beta _{n}<\frac{1}{q_{n+1}} \\
|q_{n-1}\a | _{\Z } &<&|k\a |_{\Z },~\forall ~0<k<q_{n}
\end{eqnarray*}

The following notion is essential in K.A.M. theory. It is related with the quantification of the closeness of
rational numbers to certain classes of irrational numbers.
\begin{defn} \label{def DC}
We will denote by $DC(\gamma ,\tau )$ the set of numbers $\a$ in $\T \setminus \Q$ such that for any $k\not=0$,
$|k\. \a | _{\Z }\geq \frac{\gamma ^{-1}}{|k|^{\tau }}$. Such numbers are called \textit{Diophantine}.
\end{defn}
The set $DC(\gamma ,\tau )$, for $\tau >1$ fixed and $\gamma \in \R _{+} ^{\ast}$ small is of positive
Haar measure in $\T $. If we fix $\tau$ and let $\gamma$ run through the positive real numbers, we obtain
$\cup_{\gamma >0} DC(\gamma ,\tau )$ which is of full Haar measure. The numbers that do not satisfy any Diophantine
condition are called \textit{Liouvillean}. They form a residual set of $0$ Lebesgue measure.

This last following definition concerns the relation of the approximation of an irrational number with its continued
fractions representation.
\begin{defn} \label{def RDC}
We will denote by $RDC(\gamma ,\tau )$ is the set of  \textit{recurrent Diophantine} numbers, i.e. the $\a $ in
$\T \setminus \Q$ such that $\mathrm{G}^{n}(\a)\in DC(\gamma ,\tau )$ for infinitely many $n$.
\end{defn}
The set $RDC$ is also of full measure, since the Gauss map is ergodic with respect to a smooth measure.
In contexts where the parameters $\gamma $ and $\tau $ are not significant, they will be omitted in the notation of both sets.

\section{Cocycles in $\T ^{d} \times G$} \label{sec def dynamics}

\subsection{Definition of the dynamics}

Let $\a \in \T ^{d} \equiv \R ^{d} / \Z ^{d} $, $d \in \N ^{*}$, be an irrational rotation, so that the
translation $x \mapsto x+\a \mod(\Z ^{d} ) $ is minimal and uniquely
ergodic. For the greatest part of the article $d=1$, and the results
holding for cocycles over rotations in tori of higher dimension will
be proved for $d=1$, but stated in full generality, in order to simplify notation.

If we also let $A(\. )\in C^{\infty}(\T ^{d} ,G)$, the couple
$(\a ,A(\. ))$ acts on the fibered space
$\T ^{d} \times G \ra \T ^{d} $ defining a diffeomorphism by
\begin{equation*}
(\a ,A(\. )).(x,S)=(x+\a ,A(x).S)
\end{equation*}
for any $(x,S)\in \T ^{d} \times G$. We will call such an action a
\textit{quasiperiodic cocycle over }$R_{\a}$ (or simply
a cocycle). The space of such actions is denoted by
$SW_{\a }^{\infty}(\T ^{d},G)\subset \Diff ^{\infty}(\T ^{d} \times G)$.
Most times we will abbreviate the notation to $SW_{\a }^{\infty}$.
The space of cocycles homotopic to the $\Id$ (relevant if $G$ is not
simply connected), will be denoted by $SW_{\a }^{\infty,1}$.
Cocycles are a class of fibered diffeomorphisms, since
fibers of $\T ^{d} \times G$ are mapped into fibers, and the mapping
from one fiber to another in general depends on the base point.
The number $d \in \N ^{*}$ is the number of frequencies of the cocycle.

The space $\bigcup\nolimits_{\a \in \T ^{d} }SW_{\a }^{\infty }$ will be denoted by $SW^{\infty }$.
The space $SW_{\a }^{\infty }$ inherits the topology of $C^{\infty }(\T ^{d} ,G)$, and $SW^{\infty }$ has the standard
product topology of $\T ^{d} \times C^{\infty }(\T ^{d},G)$. We note that cocycles are
defined over more general maps and in more general contexts of regularity
and structure of the basis and fibers.

If we consider a representation of $G$ on a vector space $E $, the action of the cocycle can be also defined
on $\T ^{d} \times E $, simply by replacing $S$ by a vector in $E$ and multiplication in $G $ by the action.
The particular case which will be important in this article is the representation of $G $ on $L^{2} (G)$,
and the resulting action of the cocycle on $L^{2}(\T ^{d} \times G)$.

The $n$-th iterate of the action is given by
\begin{equation*}
(\a,A(\. ))^{n}.(x,S)=(n\a,A_{n}(\. )).(x,S)=(x+n\a 
,A_{n}(x).S)
\end{equation*}%
where $A_{n}(\. )$ represents the \textit{quasiperiodic product} of matrices equal to
\begin{equation*}
A_{n}(\. )=A(\. +(n-1)\a)\cdots A(\. )
\end{equation*}%
for positive iterates. Negative iterates are found as inverses of positive
ones:%
\begin{equation*}
(\a,A(\. ))^{-n} =((\a,A(\. ))^{n})^{-1} 
=(-n\a,A^{\ast }(\. -n\a)\cdots A^{\ast }(\. -\a))
\end{equation*}

\subsection{Classes of cocycles with simple dynamics, conjugation}

The cocycle $(\a  ,A(\.  ))$ is called a constant cocycle if
$A(\. )=A\in G$ is a constant mapping. In that case,
the quasiperiodic product reduces to a simple product of matrices
\begin{equation*}
(\a  ,A)^{n}=(n\a  ,A^{n})
\end{equation*}

The group $C^{\infty }(\T ^{d},G)\hra SW^{\infty }(\T ^{d},G)$ acts by \textit{fibered conjugation}: Let
$H(\.  )\in C^{\infty }(\T ^{d},G)$ and $(\a  ,A(\. ))\in SW^{\infty }(\T ^{d},G)$. Then we define
\begin{eqnarray*}
Conj_{H(\.  )}.(\a  ,A(\.  )) &=&(0,H(\.  ))\circ (\a  ,A(\.  ))\circ (0,H(\.  ))^{-1} \\
&=& (\a  ,H(\.  +\a  ).A(\. ).H^{-1}(\.  ))
\end{eqnarray*}
which is in fact a change of variables within each fiber of the product $\T ^{d} \times G$. The dynamics of
$Conj_{H(\. )}.(\a  ,A(\. ))$ and $(\a ,A(\. ))$ are essentially the same, since
\begin{equation*}
(Conj_{H(\.  )}.(\a  ,A(\.  )))^{n}=(n\a  ,H(\.  +n\a ).A_{n}(\.  ).H^{-1}(\. ))
\end{equation*}
\begin{defn}
Two cocycles $(\a ,A(\. ))$ and $(\a ,\tilde{A}(\.  ))$ in $SW^{\infty }_{\a }$ are $H^{s}$-\textit{conjugate} iff there
exists $H(\.  )\in H^{s}(\T ^{d},G)$ such that $(\a  ,\tilde{A}(\.  ))=Conj_{H(\.  )}.(\a  ,A(\.  ))$.
We will use the notation $(\a  ,A(\.  ))\sim (\a  ,\tilde{A}(\.  ))$
to state that the two cocycles are conjugate to each other.
\end{defn}
Since constant cocycles are a class for which dynamics can be analysed, we give the following definition.
\begin{defn} \label{defn a.r.}
A cocycle will be called \textit{reducible} iff it is conjugate to a constant.
\end{defn}
Due to the fact that not all cocycles are reducible (e.g. generic cocycles in $\T \times \Sp ^{1}$
over Liouvillean rotations, but also cocycles over
Diophantine rotations, even though this result is hard to obtain, see \cite{El2002a}, \cite{Krik2001}, \cite{NKPhD})
we also need the following concept, which has proved to be crucial in the study of such dynamical systems.
\begin{defn} \label{def almost reducibility}
A cocycle $(\a ,A(\. ))$ is said to be \textit{almost reducible} if there exists a sequence of conjugations
$H_{i}(\. ) \in C^{\infty}$, such that $Conj_{H_{i}(\. )}.(\a ,A(\. ))$ becomes arbitrarily close to constants
in the $C^{\infty }$ topology, i.e. iff there exists $(A_{i})$, a sequence in $G$, such that
\begin{equation*}
A_{i}^{\ast } \left(  H_{i}(\. +\a )A(\. )H_{i}^{\ast }(\. ) \right) = e^{F_{i} (\. )} \overset{C^{\infty } }{\ra } \Id
\end{equation*}
\end{defn}
Usually, this property is established in a K.A.M. constructive way, making it possible to measure the rate of
convergence versus the explosion of the conjugations. Even though it is not a part of the definition, if
$F_{i} (\. ) \in C^{\infty} (\T ^{d} ,g )$ the error term which makes this last limit into an equality, almost reducibility
then comes along with obtaining that
\begin{equation} \label{eq rate conv}
Ad (H_{i} (\. ) ) .F_{i} (\. ) = H_{i} (\. ) .F_{i} (\. ) . H_{i}^{*} (\. ) \overset{C^{\infty } }{\ra } 0
\end{equation}
If this additional condition is satisfied, almost reducibility in the sense of
the definition above and almost reducibility in the sense that "the cocycle can be conjugated arbitrarily close
to reducible cocycles" are equivalent.

We can now recall the main theorem of \cite{NKPhD}.
\begin{thm} \label{thm global density}
Let $\a \in RDC \subset \T $ and $G $ a compact semisimple Lie group. Then, reducible cocycles are dense in
$SW_{\a }^{\infty ,1}(\T , G)$ in the $C^{\infty}$ topology.
\end{thm}
The following theorem is used the proof of the local density theorem, already proved in
\cite{KrikAst} and \cite{El2002a}.
We will sketch the proof of \cite{NKPhD}, since it implies a useful
result (lem. \ref{preliminary reduction})
in a more direct manner than the previously existing proofs.
\begin{thm} \label{thm local a.r.}
Let $\a \in DC(\gamma , \tau ) \subset \T ^{d}, d\geq 1 $ and $G$ a semisimple compact Lie group. Then, there
exists $s_{0} \in \N ^{*} $ and $\epsilon >0$, such that if
$(\a , Ae^{F(\. )} )\in SW^{\infty}_{\a} (\T ^{d}, G ) $ with $\| F(\. ) \| _{0} < \epsilon $,
$\| F(\. ) \| _{s_{0}} <1 $ and $A \in G$, $(\a , Ae^{F(\. )} )$ is almost
reducible. The neighborhood of constants thus defined will be
denoted by $\WW = \WW _{\a ,d, G}$.
\end{thm}
As a continuation of the proof of this theorem for cocycles in $\T ^{d} \times SO(3)$,
H. Eliasson, in his article cited above, proved the following one.
Originally, the theorem was proved for analytic cocycles, but the proof works equally
well in the smooth category\footnote{The genericity of H. Eliasson's condition in $C^{\infty}$ will be proved
in section \ref{sec a.r. and due}.}.
\begin{thm}[H. Eliasson] \label{thm El UE}
Let $G=SU(2)$. Then, unique ergodicity is generic in the set of cocycles satisfying the conditions
of theorem \ref{thm local a.r.}.
\end{thm}
We can also prove the following corollary.
\begin{cor} \label{global density of UE}
If $\a \in RDC$ and $G= SU(2)$, uniquely ergodic cocycles
are generic in $SW_{\a }^{\infty}(\T , G)$ in the $C^{\infty}$
topology.
\end{cor}
\begin{proof}
Every given cocycle is accumulated by reducible cocyles, by theorem
\ref{thm global density}. Reducible cocycles are accumulated by UE cocycles, by theorem
\ref{thm El UE}. Unique ergodicity is a $G_{\d }$ condition in the $C^{0 }$ topology
(this is an easy $C^{0}$ adaptation of lemma \ref{lem DUE Gd}).
\end{proof}
Finally, for the sake of completeness, we give the following theorem which belongs
to the non-local regime of the dynamics (see \cite{Krik2001},
\cite{Fra2004}, \cite{NKPhD}), which is an intermediate result in
the proof of theorem \ref{thm global density} for $G=SU(2)$.
\begin{thm} \label{normal form SU(2)}
Let $G=SU(2)$ and $\a \in RDC$. Then, all cocycles which are not
almost-reducible, are conjugate to a diagonal cocycle of the form
$(\a , E_{r}(\. ))$, where $E_{r}(\. ) = \{ e^{2i\pi r \.},0 \}$
and $r \in \N^{*}$.
\end{thm}

\section{Invariant Distributions and DUE}

\subsection{Definitions}

If we let $M$ be a compact $C^{\infty}$-manifold without boundary,
$C^{\infty}(M) =C^{\infty}(M, \C) $, the space of
$C^{\infty}$-smooth functions on $M$, is a Fr\'{e}chet space for the
$C^{\infty}$ topology. Every $f \in \Diff ^{\infty} (M)$ acts on
$C^{\infty}(M ) $ by pullback: for $\phi \in C^{\infty}(M )$ we define
\begin{equation*}
f^{*}\phi (\. ) = \phi \circ f (\. )
\end{equation*}
For a fixed $f$ this action is clearly linear and continuous. Therefore, $f$ acts by the transposed action
on the dual space of $C^{\infty}(M ) $, the space of distributions $\mathcal{D}'(M)$. This action is given by
\begin{equation*}
\langle f_{*}T , \phi \rangle = \langle T , f^{*} \phi \rangle
\end{equation*}
for all $T \in \mathcal{D}'(M)$ and for all $\phi \in C^{\infty}(M )$. We also remind that for any given distribution
$T$, there exists a minimal $k \in \N $, called the order of $T$, such that $T$ extends to a continuous functional
on $C^{k}(M )$ by density of $C^{\infty}(M ) \hra C^{k}(M )$. Distributions of order $0$ are simply the
(signed) measures of finite mass on $M$.
For any given $f \in \Diff ^{\infty} (M)$ we can define the space of $f$-invariant distributions as
\begin{equation*}
\mathcal{D}^{\prime}_{f} (M) = \{ T\in \mathcal{D}'(M) , f_{*}T = T \}
\end{equation*}
This space is non-empty, since it contains the (vector space generated by the) invariant measures of $f$, so it
has dimension at least $1$.
\begin{defn}
A diffeomorphism $f$ is DUE iff $\mathcal{D}^{\prime}_{f} (M)$ is one-dimensional.
\end{defn}
In particular, Distributional Unique Ergodicity implies Unique Ergodicity. The property of being DUE, however, refines
UE, since a space much larger than the space of Radon measures, that of distributions, is considered.

We will use an alternative definition of DUE, stated in terms of the cohomological equation. It is a direct
consequence of the Hahn-Banach theorem. The space
$C^{\infty }_{0} (M) $ is the space of $0$-mean-value functions on $M$ for the unique ergodic measure of $f$.
\begin{lem} \label{lem DUE by approx}
A diffeomorphism $f \in \Diff ^{\infty} (M) $ is DUE iff for every $\phi (\. ) \in C^{\infty }_{0} (M) $ and every
$\d ,s_{0} > 0$, there exist $\e (\. ) \in C^{\infty } (M) $,
satisfying $\| \e (\. ) \| _{s_{0}} \leq \d $, and
$\psi \in C^{\infty } (M)$ such that
\begin{equation*}
\phi (\. ) = f^{*} \psi (\. ) - \psi (\. ) + \e (\. )
\end{equation*}
\end{lem}
We give the following definition, giving a name to the functions appearing on the rhs of the previous
equation.
\begin{defn}
Functions of the form
\begin{equation*}
\{ f^{*} \psi (\. ) - \psi (\. ) , \psi \in C^{\infty }(M,\C )\}
\end{equation*}
are called coboundaries and are denoted by $Cob^{\infty }(f)$.
\end{defn}
Therefore, the lemma states that, $f$ is DUE iff coboundaries
are dense in codimension one subspace of $C^{\infty }(M)$. The proof can be found in \cite{AFKo2015}. Naturally,
application of the Hahn-Banach theorem cannot infer anything about whether the space of coboundaries
is closed or not. This is an additional property of the diffeomorphism $f$, for which we give
the following definition.
\begin{defn}
If the space of coboundaries is
closed, then $f $ is called \textit{cohomologically stable}.
\end{defn}

It is a
classical fact that Diophantine rotations in
tori are both DUE and cohomologically stable, while Liouvillean
rotations are only DUE. The distinction is captured by the following definition.
\begin{defn}
A diffeomorphism is called Cohomologically Rigid iff it is both DUE and cohomologically
stable.
\end{defn}

The following conjecture can be found in \cite{HurdKatConj} (problem 17).
\begin{conj} [A. Katok] \label{conj Katok}
If $f \in \Diff  ^{\infty} (M) $ is DUE and cohomologically stable, then it is smoothly conjugate to a
Diophantine rotation in a torus $\T ^{d }$.
\end{conj}
This conjecture has only been proved partially, see \cite{Koc09},
\cite{Koc07} and references therein. In the present article, we establish the
validity of the conjecture in the space of quasi-periodic cocycles in $\T \times SU(2)$
over an $RDC $ rotation.

We note that it is possible that $f$ be uniquely ergodic, but
$\dim (\mathcal{D}_{f}'(M) ) >1$, see \cite{Kat01} for
the study of the parabolic map
\begin{equation*}
\begin{array}{ccc}
\T \times \Sp ^{1} & \ra & \T \times \Sp ^{1} \\
(x ,e^{2i \pi \theta}) & \mapsto  & (x + \a, e^{2i \pi ( \theta +x )})
\end{array}
\end{equation*}
This map corresponds, in our context, to the models of not almost reducible
dynamics, $(\a ,E_{r}(\. ))$ (cf. theorem \ref{normal form SU(2)}).

\subsection{DUE cocycles}

The part of the theorem of \cite{AFKo2015} which we are interested in reads as follows.
\begin{thm} [A.Avila-B.Fayad-A.Kocsard] \label{DUE AFKo}
DUE diffeomorphisms are generic in $\mathcal{AK} ^{\infty} \subset \Diff ^{\infty } (\T \times P) $
defined as the closure of 
\begin{equation}
\{ Conj_{B(\.) } (\a , \Id ) , \a \in \T , B(\. ) \in C^{\infty} (\T , P )  \}
\end{equation}
in $\Diff ^{\infty} (\T \times P)$. Here, $P = G / H $ with $G$ a compact Lie group and $H$ a closed subgroup.
\end{thm}
We will also use the following fact, proved in the same reference (as well as in
\cite{NKCohomRigAK}).
\begin{lem} \label{lem DUE Gd}
The set of DUE diffeomorphisms of a compact manifold $M$
is a $G_{\d }$ subset of $\Diff ^{\infty}(M)$.
\end{lem}

%

\subsection{Establishing DUE in the quasi-periodic setting}

The definition of DUE being difficult to work with, one can adopt a more constructive
approach to tackling the problem, especially in the quasi-periodic setting, or in the
presence of algebraic structure more generally. Our strategy for studying invariant
distributions of cocycles is encapsulated in the following proposition.

\begin{prop} \label{prop DUE qp}
Let $f \in \Diff  ^{\infty } (M)$ and suppose that there exists a sequence
$\{ f_{n} \} \subset \Diff  ^{\infty } (M) $, $f_{n} \ra f$, satisfying the following properties.
For every function $\f \in \smm$, for every
$s \in \N $ and for every $\epsilon >0$, there
exists $n \in \N$ and a function $\f _{n} \in \smm$ such that
\begin{enumerate}
\item $\| \f - \f _{n} \| _{s} < \epsilon $
\item There exists $\psi _{n} \in \smm$ such that
\begin{equation} \label{eq DUE approx}
\psi _{n} \circ f_{n} - \psi_{n} = \f _{n}
\end{equation}
(i.e. $\f _{n} \in Cob^{\infty}(f_{n})$)
\item The solution $\psi _{n}$ satisfies the following estimate
\begin{equation} \label{eq error terms}
C_{s }\| \psi _{n} \| _{s+1} \vertiii{f}_{s} d_{s}(f,f_{n}) < \epsilon
\end{equation}
\end{enumerate}
Then $f$ is $DUE$.
\end{prop}
%

\begin{proof}
The proposition is practically a tautology, since we can rewrite eq. \ref{eq DUE approx} as
\begeq
\psi _{n} \circ f - \psi _{n} = \f + (\f _{n} - \f ) + ( \psi _{n} \circ f - \psi _{n}\circ f_{n})
\endeq
Then, the estimate of eq. \ref{composition} implies that the terms in the parentheses
sum up to $<2 \e $ in the $H^{s}$ norm by assumption. This is equivalent to $f$ being
DUE by lemma \ref{lem DUE by approx}.
\end{proof}

This approach implicitly assumes sufficiently good understanding of the coboundary
operator for the approximants $f_{n}$. In the quasiperiodic setting, this is
achieved by approximating the given system $f$ by systems $f_{n}$ that are conjugated to
a simple system, in our case a constant (moreover, a resonant -cf. the proof of lemma
\ref{loc conj lem}) cocycle.

The cohomological equation transforms naturally under changes of coordinates, thanks to the
following obvious lemma.
\begin{lem} \label{lem cohom eq change coord}
If $f: M\ra M$ is
a $\sm$-smooth diffeomorphism and
\begeq
\f = \psi \circ f - \psi
\endeq
then, if $\tilde{f}= g \circ f \circ g^{-1}$, we obviously have
\begeq
\f \circ g^{-1} = \psi \circ g^{-1} \circ \tilde{f} - \psi \circ g^{-1}
\endeq
\end{lem}
This property could prove itself useful for showing that a diffeomorphism is Cohomologically
Rigid. However, in our case we will establish that our examples are DUE and not rigid, so we
need to look into what happens when we transform approximate solutions.

If for fixed $f_{n}$ and a given function $\f$ we can obtain an approximate solution to the
cohomological equation,
\begeq
\f = \psi \circ f_{n} - \psi + \e
\endeq
with $\e $ a small $\sm$ function, then when we change coordinates as above, we
obtain
\begeq
\f \circ g^{-1} = \psi \circ g^{-1} \circ \tilde{f}_{n} - \psi \circ g^{-1} + \e
\circ g^{-1}
\endeq
This approximate solution will give useful information only if we can control the size
of $\e \circ g^{-1}$, which in view of eq. \ref{composition} boils down to limiting the
size of $g$ with respect to $\e$. This, alongside with the fact that
in out case the
coboundary space of the approximats is of infinite codimension,
is the main subtlety of our approximation argument laid out in \S
\ref{section proof main thm}.


For a deeper study of this approach to establishing DUE in the quasi-periodic
setting we refer the reader to \cite{NKCohomRigAK}.

\subsection{Concerning the conjecture by A. Katok}

Theorem \ref{thm DUE} does not provide counterexamples to conjecture \ref{conj Katok},
because of cor. \ref{no cout-ex}.
The heuristics for this are the following. Cocycles can be DUE
if the dynamics in the fibers look
like a Liouvillean (or resonant) rotation around a fixed axis for
some time long
enough for equidistribution of orbits in each leaf
$\approx \T \times \Sp ^{1}$ to be established to a certain extent.
Then, the axis tilts, and the procedure continues for an infinite
number of times. The Liouvillean character of the dynamics rules out
cohomological stability, and an additional generic condition on the
angles between the successive constants, stricter than the one needed
for unique ergodicity, is needed in order to establish DUE.

\subsection{A global theorem} \label{some direct cons}
The proof of theorem \ref{DUE AFKo}, combined with theorem \ref{thm global density}, gives
the following corollary
which identifies the space where DUE lives, showing that DUE is not an exotic property in
these non-abelian extensions of rotations. We recall the notation $SW ^{\infty ,1}(\T , G)$
for cocycles that are homotopic to the $\Id$.
\begin{cor} \label{genericity in hom to csts}
DUE is a $G_{\delta}$-dense property in $SW ^{\infty ,1}(\T , G)$,
where $G$ is a semisimple compact Lie group.
\end{cor}
This is the only part of the article where we consider Lie groups
other than $SU(2)$, and a reader unfamiliar with the theory of
semi-simple compact Lie groups can replace $G$ with $SU(N)$.
\begin{proof}
We begin by observing that theorem \ref{thm global density}, together with the density of $RDC$ in $\T $,
implies that reducible cocycles are dense in $SW ^{\infty ,1} (\T , G)$, since we allow the frequency to vary.
Secondly, we point out that the proof of theorem \ref{DUE AFKo} actually shows that DUE is
$G_{\delta}$-dense in the Anosov-Katok space
\begin{equation*}
\mathcal{AK} ^{\infty} =  \{ Conj_{B(\.) } (\a , \Id ) , \a \in RDC , B(\. ) \in C^{\infty} (\T , G )  \}^{cl_{\infty}}
\end{equation*}
The only part missing in the proof is actually showing that this space coincides with the space of
cocycles with values in $G$, which are homotopic to the $\Id$.

In what follows, we drop the rigor in order to skip some technicalities on semisimple compact Lie groups. Let us
fix an embedding $\T ^{w} \equiv \TT \subset G $ such that $w \in \N ^{*} $ is maximal. Therefore, we also fix
$\mathcal{Z}$, the lattice of preimages of $\Id \in G $ under the $\exp$
in $\R ^{w} \equiv T_{0} \T ^{w}$. Then, we see that
reducibility of a cocycle to $(\a ,\Id )$ with $\a \in \T \setminus \Q $ implies reducibility to $(\a , \exp (s))$, where
$s \in \R ^{w} $ satisfies $s \in  \a \mathcal{Z} $. This can be seen by considering
\begin{equation*}
\{ Conj_{B(\.)}(\a ,\Id ) , B(\. ) : \T \ra \TT \text{, a morphism} \}
\end{equation*}
By construction, $B(\a )$ satisfies the hypothesis above. Finally, fixing a torus
$\TT $ does not affect generality, since all constant cocycles live in the fixed torus, modulo conjugation by a constant
(\cite{DieudElV}, \cite{DuisterKolkLieGroups}). Since, by irrationality of $\a $, the
exponentials of such $s\in  \a \mathcal{Z}$ are dense in $\TT $, we have
\begin{equation*}
\mathcal{AK} ^{\infty} =  \{ \cup _{\a \in RDC} \{
(\a ,B(\. +\a ).A.B^{*}(\. ), A \in G, B(\. ) \in C^{\infty}
(\T ,G)) \}  \}^{cl_{\infty}}
\end{equation*}
Applying theorem \ref{thm global density} for each fixed $\a \in RDC $ gives that
\begin{equation*}
\mathcal{AK} ^{\infty} =  \{ \cup _{\a \in RDC} SW_{\a }^{\infty ,1}(\T , G)  \}^{cl_{\infty}}
\end{equation*}
This last space by our first remark is $SW ^{\infty ,1}(\T , G)$.
\end{proof}
This corollary shows that the local picture as obtained in theorem \ref{thm El UE} is actually the corresponding part of the
global picture, where uniquely ergodic, and therefore not reducible, cocycles form a generic subset in the space of cocycles.
We remind, nonetheless, that we have allowed the frequency to vary
(compare with theorem \ref{thm El UE}, which holds for a fixed
frequency, and in fact for an arbitrary number of frequencies). The
analogy is, consequently, quite loose.

\section{Review of almost reducibility and unique ergodicity} \label{sec a.r. and due}

In this section we present the basic points of the proof of theorem \ref{thm local a.r.},
which is the basis for proving theorems \ref{thm El UE} and \ref{finite reduc thm}.
In fact, we will not (re-)prove theorem \ref{thm El UE}, since our goal is to prove a stronger one\footnote{Though in the
$C^{\infty } $ category, and not in the real analytic one.}.

The rest of this section is devoted to proving theorem \ref{thm local a.r.}
and analyzing the conjugations produced by the K.A.M. scheme used in
its proof. From now on,
\begeq
G = SU(2)\text{ and } d=1
\endeq
even though the local
almost reducibility theorem (and therefore also its direct consequences
not requiring renormalization of the dynamics in their proof) is true
for an arbitrary number of frequencies $d \in \N^{*}$.

The proof of theorem \ref{thm local a.r.} is by K.A.M. iteration of lemma
\ref{loc conj lem}, and the proof of theorem \ref{finite reduc thm} is based on a
more careful analysis of this K.A.M. iteration.

The notation introduced here will be used throughout the rest of the paper.

\subsection{Local conjugation}
Let $(\a , Ae^{F(\. )} ) =(\a , A_{1}e^{F_{1}(\. )} ) \in SW^{\infty} (\T ,G) $ be a cocycle over a
Diophantine rotation satisfying some smallness conditions to be made more precise rightaway. Without any loss of
generality, we can also suppose that $A = \{ e^{2i\pi a} , 0 \} $ is diagonal. The goal is to conjugate the cocycle
ever closer to constant cocyles by means of an iterative scheme. This is obtained by iterating the following lemma,
for the detailed proof of which we refer to \cite{KrikAst}, \cite{El2002a} or \cite{NKPhD}. The rest of this
paragraph is devoted to a summary of the proof, for the sake of completeness.
\begin{lem} \label{loc conj lem}
Let $\a  \in DC(\gamma ,\tau )$ and $K\geq C \gamma N^{\tau +1/2 }$. Let, also,
$(\a  ,Ae^{F_{1}(\. )})\in SW^{\infty }(\T ,G)$ with
\begin{equation*}
c_{1}KN ^{s_{0}} \e  _{1,0}<1
\end{equation*}
where $c_{0}, s_{0} $ depend on $\gamma , \tau $ (and $d$), and
$\e  _{1,s}=\left\Vert F_{1}\right\Vert _{s}$. Then, there exists a
conjugation $G_{1}(\. )=G(\. )\in C^{\infty }(\T ,G)$ such that
\begin{equation*}
G(\.  +\a  ).A_{1}.e^{F_{1}(\.  )}.G^{*}(\.  )=A_{2} e^{F_{2} (\.  )}
\end{equation*}
and such that the mappings $G(\. )$ and $F_{2}(\. )$ satisfy the
following estimates
\begin{eqnarray*}
\left\Vert G(\.  )\right\Vert _{s} &\leq &  c_{1,s} (N^{s} + KN^{s+1/2}\e  _{1,0}) \\
\e  _{2,s}&\leq &c_{2,s}K^{2}N^{2\tau +1}(N^{s}\e  _{1,0}+\e _{1,s})\e  _{1,0}+
C_{s,s^{\prime }}K^{2}N^{s-s^{\prime }+2\tau +1}\e _{1,s^{\prime }}
\end{eqnarray*}
where $s' \geq s$.
\end{lem}

We now sketch the proof of the lemma.

If we suppose that $Y(\. ) : \T \ra g $ can conjugate $(\a , A_{1}e^{F_{1}(\. )} )$ to
$(\a , A_{2}e^{F_{2}(\. )} )$, with $\| F_{2}(\. ) \| \ll \| F_{1}(\. ) \| $, then it must satisfy the functional
equation
\begeq
\begin{array}{r@{}l}
e^{Y(\. +\a )}A_{1}e^{F_{1}(\. )}e^{-Y(\. )}  &= A_{2}e^{F_{2}(\. )}
\Leftrightarrow
\\
A_{1}^{*}e^{Y(\. +\a )}A_{1}e^{F_{1}(\. )}e^{-Y(\. )}  &= A_{1}^{*} A_{2}e^{F_{2}(\. )} 
\end{array}
\endeq
Linearization of this equation under the assumption that all $C^{0}$ norms are smaller than $1$ gives
\begin{equation*} 
Ad(A_{1}^{*})Y(\. +\a ) +F_{1}(\. )-Y(\. )  = \exp ^{-1} (A_{1}^{*} A_{2})
\end{equation*}
which we will write in the eigenspaces of $Ad(A_{1}^{*})$,
separating the diagonal from the non-diagonal part.

The equation for the diagonal coordinate reads
\begin{equation*}
Y_{t } (\. +\a )-Y_{t } (\.) +F_{1,t }(\. )=0
\end{equation*}
This equation can be solved by considering Fourier coefficients. For reasons well known in
K.A.M. theory, we
truncate at an order $N \in \N ^{*} $ to be determined by the parameters of the problem and obtain a solution
to the equation
\begin{equation*}
Y_{t } (\. +\a )-Y_{t } (\.)+\dot{T}_{N} F_{1,t }(\. )=0
\end{equation*}
satisfying the estimate
\begin{equation*}\| Y_{t } (\.) \| _{s} \leq \gamma C_{s} N^{s + \tau + 1/2} \e _{1,0}
\end{equation*}
The rest satisfies the estimate of eq. \ref{truncation est}. The mean value $\hat{F}_{1,t } (0) $ is an obstruction and
will be integrated in $\exp ^{-1} (A_{1}^{*} A_{2}) $.

As for the equation concerning the non-diagonal part, it reads
\begin{equation} \label{untruncated local eq}
e^{-4i\pi a}Y_{z }(\. +\a )-Y_{z }(\. )+F_{1,z}(\. )=0
\end{equation}
or, in the frequency domain,
\begin{equation} \label{local eq in Fourier}
(e^{2i\pi ( k\a  - 2a )}-1)\hat{Y}_{z  }(k)=-\hat{F}_{1,z }(k),~ k \in \Z
\end{equation}
Therefore, the Fourier coefficient $\hat{F}_{1,z }(k_{r })$ cannot be eliminated with good estimates if
\begin{equation*}
| k_{r } \a  -2a | _{\Z } < K ^{-1}
\end{equation*}
for some $K>0 $ big enough. If $K = N^{\nu }$, with $\nu > \tau $,
then we know by \cite{El2002a} (see also lemma \ref{Dioph csts}) that
such a $k_{r }$ (called a \textit{resonant mode}), if it exists and
satisfies $0\leq k_{r} \leq N $, is unique in
$\{ k \in \Z , |k- k_{r } | \leq 2N \}$. In other words, for
$0 < |k- k_{r } | \leq 2N$,
\begeq
| k_{r } \a  -2a | _{\Z } \geq K ^{-1} = N^{-\nu }
\endeq

We now suppose that such a $k_{r}$ exists and is non-$0$, and
call $T^{k_{r }} _{2N} $ the truncation operator
projecting on the frequencies $ 0 < |k- k_{r } | \leq 2N$.
Thus the equation
\begin{equation*}
e^{-4i\pi a}Y_{z }(\. +\a )-Y_{z }(\. )=-T^{k_{r }} _{2N} F_{1,z}(\. )
\end{equation*}
can be solved and the solution satisfies
\begeq
\| Y_{z } (\.) \| _{s} \leq C_{s} N^{s + \nu + 1/2} \e _{1,0}
\endeq
We will define the rest operator $R^{k_{r }} _{2N}$ in the obvious way, projecting on
frequencies $|k- k_{r } | > 2N$. The spectral support of $  F_{1,z} (\. )$ is thus
outside $[-N,N]\cap \Z $, and can therefore be estimated like a classical rest operator
$R_{N}$. If such a $k_{r }$ does not exist, we project on $ |k| \leq 2N$, but this case
more straightforward and left to the reader.

In total, the equation that can be solved with good estimates is
\begin{equation} \label{eq solved lin eq}
Ad(A_{1}^{*})Y(\. +\a ) -Y(\. ) + F_{1}(\. )= 
\{ \hat{F}_{1,t } (0), \hat{F} _{1,z} (k_{r }) e^{2i\pi k_{r } \.} \} +
\{ R_{N} F_{1,t} (\. ) , R^{k_{r }} _{2N} F_{1,z} (\. ) \}
\end{equation}
with $\| Y(\. ) \|_{s} \leq C_{s} N^{s + \nu + 1/2 } \e _{1,0} $.
Under the smallness assumptions of the hypothesis, the linearization
error is small and the conjugation thus constructed satisfies
\begin{equation} \label{eq solved non-lin eq}
\begin{array}{r@{}l}
e^{Y(\. +\a )}A_{1}e^{F_{1}(\. )}e^{-Y(\. )} 
&=
\{ e^{2i\pi (a + \hat{F}_{1,t } (0) )} ,0 \} _{G} . e^{\{ 0, \hat{F} _{1,z} (k_{r })
e^{2i\pi k_{r } \.} \} _{g} } e^{F_{2}^{\prime}(\. )}
\\
&=
\{ e^{2i\pi (a + \hat{F}_{1,t } (0) )} ,0 \} _{G} . e^{\tilde{F} _{1} (\.)} e^{F_{2}^{\prime}(\. )}
\end{array}
\end{equation}
with $F_{2 }^{\prime} (\. ) $ a quadratic term. We remark that, a priori, the obstruction
$\tilde{F} _{1} (\.)$ is of the order of the initial perturbation and therefore
what we called $\exp ^{-1} (A_{1}^{*} A_{2}) $ is actually not constant in the presence of
resonant modes.

If $k_{r }$ exists and is non-zero, iteration of the lemma as it has
been proved so far is impossible. On the other hand, the conjugation
$B(\. ) = \{ e^{- 2i\pi k_{r } \. /2} ,0 \} $ satisfies
\begin{equation} \label{eq red res mode}
\begin{array}{r@{}l}
Conj_{B(\. )} (\a , A_{1}.\exp (\{ \hat{F}_{t } (0), \hat{F} _{z} (k_{r }) e^{ 2i\pi k_{r } \.}) \} ) &=
(\a , A_{1}^{\prime}.\exp (\{ \hat{F}_{t } (0), \hat{F} _{z} (k_{r })  \} ) \\
 &= (\a , A_{2})
\end{array}
\end{equation}
where we have called
$A_{1}^{\prime} = B(\a ) A_{1} = \{ e^{2i\pi (a - k_{r} \a /2)} ,0\}$.
That is, $B(\. )$ reduces the initial constant perturbed by the obstructions to a cocycle close to $(\a , \pm \Id )$.

Application of the conjugation $B(\. )$ to eq. \ref{eq solved lin eq} gives
\begin{eqnarray*}
Ad((A'_{1})^{*})Y'(\. +\a ) -Y'(\. ) + F_{1}^{\prime}(\. ) &=& 
\{ \hat{F}_{t } (0), \hat{F} _{z} (k_{r })  \} +
\{ R_{N} F_{t} (\. ) , e^{- 2i\pi k_{r } \.} R^{k_{r }} _{2N} F_{z} (\. ) \} \\
&=& \{ \hat{F}_{1,t }^{\prime} (0), \hat{F} _{1,z}^{\prime} (0)  \} +
\{ R_{N} F_{1,t}^{\prime} (\. ) , R _{2N} F'_{1,z} (\. ) \}
\end{eqnarray*}
where we called
$F_{1}^{\prime} (\. ) = Ad(B(\. )).F_{1}(\. ) = \{ F_{1,t}(\. ) , e^{-2i\pi _{k_{r} \.} } F_{1,z} (\. ) \} $, similarly for $Y(\. )$.

The equation for primed variables can be
obtained from eq. \ref{untruncated local eq} by applying $Ad(B(\. ))$ and using the fact
that $B(\. )$ is a
morphism commuting with $A_{1}$.

Applying the same conjugation $B(\. )$ to the non-linear problem, eq.
\ref{eq solved non-lin eq}, gives
\begin{equation} \label{eq solved conj pb}
\begin{array}{r@{}l}
Conj_{B(\. ).e^{Y(\. )}}A_{1}e^{F_{1}(\. )}&=
A_{1}'
.e^{{\{\hat{F}_{1,t } (0) }, \hat{F} _{1,z} (k_{r }) \} _{g} }
e^{F_{2}(\. )}\\
&= A_{2}.e^{F_{2}(\. )}
\end{array}
\end{equation}
where $F_{2}(\. )$ satisfies the estimates of the lemma. We can conclude the proof by
calling
\begeq
G(\. ) = B(\. ) e^{Y(\. )}
\endeq

There is a slight complication, as $B(\. )$ may be $2$-periodic. If it is so, we can conjugate a second time with a
minimal geodesic $C(\. ): 2 \T \ra G$ commuting with $A_{1}'$ and such that $C(1) = -\Id $. The cocycle that we obtain in this
way is close to $(\a , \{ e^{i\p \a} ,0 \} _{G} )$ and the conjugation
is $1$-periodic.

We will omit reference to the conjugation $C(\. )$ from now on, since it
does not affect the estimates. In \cite{NKContSpec} we showed that this is just an
artefact of the construction.

\subsection{K.A.M. iteration and convergence}
Lemma \ref{loc conj lem} can serve as the step of a K.A.M. scheme, with the following
standard choice of parameters:
\begin{equation} \label{eq def Nn}
N_{n+1} = N_{n}^{1+\s } = N^{(1+\s )^{n-1}}
\end{equation}
where $N=N_{1}$ is big enough and
$0<\s <1$, and $K_{n} = N_{n}^{\nu }$, for some $\nu > \tau$. If we suppose that $(\a ,A_{n}e^{F_{n}(\. )})$
satisfies the hypotheses of lemma \ref{loc conj lem} for the corresponding parameters, then we obtain
a mapping $G_{n}(\. ) = B_{n}(\. ) e^{Y_{n}(\. )} $ that conjugates it to $(\a ,A_{n+1}e^{F_{n+1}(\. )})$,
and we use the notation $\e _{n,s} = \| F_{n} \|_{s}$.

Assuming that the initial perturbation is small in small norm: $\e_{1,0} < \epsilon <1$, and not big in
some bigger norm: $\e _{1, s_{0}}<1$, where $\epsilon$ and $s_{0}$ depend on the choice
of parameters,\footnote{The smallness condition is irrelevant in the present article,
and the letters $\epsilon$ and $s_{0}$ will be used with a different meaning.}
we can prove (see \cite{NKPhD} and, through it, \cite{FK2009}), that the scheme can be iterated,
and moreover
\begin{eqnarray*}
\e _{n,s} &=& O(N_{n}^{-\infty}) \text{ for every fixed } s \text{ and }\\
\| G_{n} \| _{s} &=& O (N_{n}^{s+\lambda}) \text{ for every } s \text{ and some fixed } \lambda >0
\end{eqnarray*}
We say that the norms of perturbations decay exponentially, i.e.
\begeq
N_{n}^{l}\e _{n,s} \ra 0 , \forall s,l \in \N
\endeq
while conjugations grow at most polynomially. This proves theorem \ref{thm local a.r.}.

\subsection{The K.A.M. normal form} \label{remarks and cons}
The product of conjugations $H_{n} = G_{n} \cdots G_{1} $, which by construction satisfies
$Conj_{H_{n}} (\a ,A_{1} e^{F_{1}}) = (\a ,A_{n+1} e^{F_{n+1}})$, is not expected to converge. In fact,
it converges iff $B_{n} (\. ) \equiv \Id $, but for a finite number of steps. We thus
have the following obvious lemma.
\begin{lem} \label{lemconv conj smooth}
The sequence of conjugations constructed by the KAM scheme, $H_{n}(\. )$ converges if,
and only if, it converges in the $\sm $ topology.
\end{lem}

We refer to \cite{NKRigidity}
for the proof of the convergence of $H_{n}$ if the cocycle is
$L^{2}$-conjugate to a Diophantine constant cocycle; and of the
generic divergence of $H_{n}$ if the constant is Liouville.
Therefore, reducibility cannot be ruled out by the
divergence of $H_{n}$, at least without the further analysis of the
scheme carried out in the following section, \S
\ref{loc dens of fin dif} leading to the proof of theorem \ref{finite reduc thm}.
Anyhow, we can obtain a \textit{K.A.M. normal form} for cocycles close to constants,
which will greatly facilitate the proof of theorem \ref{finite reduc thm}.

\begin{lem} \label{preliminary reduction}
Under the hypotheses of theorem \ref{thm local a.r.},
there exists $D(\. ) \in C^{\infty } (\T ,G)$
such that, if we call $Conj_{D(\. )} (\a ,A e^{F}) = (\a ,A' e^{F'})$, then the K.A.M. scheme applied to
$(\a ,A' e^{F'})$ for the same choice of parameters consists only in the reduction of resonant modes: The
resulting product of conjugations $H_{n_{i}}(\. )$ has the form
$\prod \nolimits _{i} ^{1} B_{n_{j}} (\. ) $, where
$\{ n_{i} \} $ are the steps in which reduction of a resonant mode took place.
\end{lem}

All cocycles henceforth will be considered to be given directly in K.A.M. normal form,
and reference to the conjugation $D(\. )$ will be omitted.

In fact, we implicitly suppose an intelligent scheme, where, if we separate the perturbation $e^{F_{n}(\.)}$  in
$e^{\tilde{F}_{n}(\.)}.e^{F_{n,nr}(\.)}$, the resonant and non-resonant part, and the non-resonant part
satisfies the inductive smallness conditions of step $n+1$, the scheme
produces only the reduction of the resonant part.
Of course, nothing changes as to the significant part of the scheme.
\begin{proof}
The product of conjugations produced by the scheme at the $n$-th step is written in the form
$H_{n}(\.  ) = B_{n}(\.  )e^{Y_{n}(\.  )}...B_{1}(\. )e^{Y_{1}(\.  )}$, where the $B_{j}(\.  )$ reduce the resonant modes or $=\Id$.
We can rewrite the product in the form
\begin{equation*}
B_{n}(\.  )\dots B_{1}(\. ).e^{\tilde{Y}_{n}(\.  )} \cdots e^{\tilde{Y}_{2}(\.  )} e^{Y_{1}(\.  )}
\end{equation*}
where $\tilde{Y}_{j}(\.  ) = \prod _{j-1}^{1} Ad(B_{i}^{*}(\. )).  Y_{j}(\.  )$. Since the $Y_{j}(\.  )$ converge
exponentially fast to $0$ (they are conjugations comparable with $F_{j}$ with a fixed loss of derivatives)
in $C^{\infty}$, and
since the algebraic conjugation deteriorates the $C^{s}$ norms by a factor of the order of $N_{n-1}^{s+1}$,
$\prod _{\infty}^{1} \exp( \tilde{Y}_{j}(\.  ))$ always converges, even if the $H_{n}(\. )$ do not.
The rest of the properties follow by construction.
\end{proof}

\begin{notation}
For a cocycle in normal form, we relabel the indexes as
\begin{equation*}
(\a ,A_{n_{i}}e^{F_{n_{i}}}) = (\a ,A_{i}e^{F_{i}})
\end{equation*}
\end{notation}

The lemma therefore asserts that any cocycle in $\WW $
can be conjugated to a cocycle for which the following holds.
At every step of the K.A.M. scheme, the image of the mapping
\begeq
x \mapsto A_{i}e^{F_{i}(x)}
\endeq
is such that the argument of the eigenvalues, $a_{i}$, is close to some $k_{i} \a /2 $, and therefore not $0$,
so that a direction in $g$ is defined by it. To the first order, the image of the path $A_{i}e^{F_{i}}$ is a circle centered
on $A_{i}$ and contained in the plane orthogonal to the direction of $a_{i}$, and it is traversed $k_{i}$ times
when $x$ goes once around $\T $ (see also the figure in \cite{NKContSpec}).

The significant parameters in this expression are:
\begin{enumerate}
\item the resonant mode $k_{i} \in \Z ^{d}$. It is of the order of $N_{i}$
\item the distance of the eigenvalues from the exact resonance $\epsilon _{i} = a_{i} - k_{i} \a /2 \mod \Z$. It is at most of the order of $N_{i}^{-\nu}$
\item the size of the resonant mode $| \hat{F}_{i} (k_{i}) | \neq 0$. It is of the order
of $O(N_{i}^{-\infty})$
\end{enumerate}
The argument of the complex number $\hat{F}_{i} (k_{i})$
is not important, since, if $D\in G $ commutes with $A_{i}$, $Conj_{D}(\a ,A_{i}e^{F_{i}}) = (\a ,A_{i}e^{Ad(D).F_{i}})$ and
the perturbation term has a resonant Fourier coefficient $k_{i}$
with the same size, but a different argument.

The change of coordinates
$(x,S) \mapsto (x, B_{i}(x).S) $ sends this circle to the constant:
\begin{eqnarray*}
A_{i+1} &=& \{ e^{2i\pi \epsilon _{i}} , 0 \}_{G} . e^{\{ 0 , \hat{F}_{i}(k_{i}) \}_{g}} \\
&=& e^{ \{ 2\pi \epsilon _{i} , \hat{F}_{i}(k_{i}) \}_{g} +
O(\epsilon _{i}.|\hat{F}_{i}(k_{i})|)}
\end{eqnarray*}
in coordinates where $A_{i}$ is a diagonal matrix.
This is the meaning of eq. \ref{eq red res mode}.

Using these parameters, we can estimate the non-commutativity between
$A_{i}$ and $A_{i+1}$ using the following lemma.
\begin{lem} \label{lem est diag}
Let $(A_{i})$ be the sequence of constants in $G$ produced by the KAM scheme built upon
lemma \ref{loc conj lem}. Then, if $D_{i} \in G$ diagonalizes
$A_{i+1}$ in a basis where $A_{i}$ is diagonal, then
\begeq
d(D_{i},\Id) \simeq  \sin (\theta _{i}/2 )
\endeq
where $\theta _{i}$ is the angle formed by preimages of $A_{i}$ and $A_{i+1}$ in $su(2)$
under the exponential.
\end{lem}
\begin{proof}
Chose coordinates where $A_{i}$ is diagonal.
Direct calculation shows that, if we call $\RR (\theta ) = \{ \cos \theta , i \sin \theta \}_{SU(2)}$,
then
\begeq
\begin{array}{l@{}l}
\RR (\theta ). \{ i , 0 \}_{su(2)} . \RR ^{*}(\theta ) &= \{ i \cos 2\theta , \sin 2 \theta \}_{su(2)}
\\
\RR (\theta ). \{ 0 , i \}_{su(2)} . \RR ^{*}(\theta ) &= \{ 0, i \}_{su(2)}
\end{array}
\endeq
and that, if we call $\EE (\f ) = \{ e^{i\f },0 \}_{SU(2)}$, then
\begeq
\EE (\f ). \{ 0, 1 \}_{su(2)} . \EE ^{*}(\f ) = \{ 0 , e^{2i \f } \}_{su(2)}
\endeq
Consequently, if $\arg \hat{F}_{i}(k_{i}) = \f _{i}$ and if
\begin{equation*}
\arctan \frac{|\hat{F}_{i}(k_{i})|}{\epsilon _{i}} = \theta _{i}
\end{equation*}
the conjugation
\begeq
\begin{array}{r@{}l}
D_{i}&= \EE (\f _{i}/2) \RR (-\theta _{i}/2 ) \EE (- \f _{i}/2)\\
&= \{\cos \theta _{i} /2 , i e^{2i\pi \f _{i}} \sin \theta _{i} /2 \}_{SU(2)}
\end{array}
\endeq
diagonalizes
$A_{i+1}$ in a basis where $A_{i}$ is diagonal and, moreover, its distance to the $\Id$ is
minimal among such matrices, namely
\begeq
d(D_{i},\NN _{i}) = |\sin \theta _{i}/2|
\endeq
This concludes the proof of the lemma
\end{proof}

As a consequence,
\begeq
d_{G } ( A_{i+1} ,\NN _{i} ) \simeq \sin ( \theta _{i} /2 )
\endeq
where $\NN_{i}$ is the normalizer of the unique maximal torus passing by $A_{i}$.
The hypothesis that the cocycle is in normal form implies that $a_{i+1} $ is
$N_{n_{i+1}}^{-\nu }$-close to some resonant sphere of radius
$|k_{i+1}\a /2 |_{\frac{1}{2} \Z }$ and that the same description holds at the next step,
unless the K.A.M. scheme stops reducing resonances after a finite number of steps. In that
case, the K.A.M. form is finite and the cocycle is $\sm$-reducible.

\subsection{Unique ergodicity}

This K.A.M. normal form allows us to discuss the genericity of UE for $C^{\infty}$-smooth
cocycles, theorem \ref{thm El UE}, in the notation we have established.
Let us introduce the notation
\begeq
LA(\. ) = \partial A(\. ) A^{*}(\. ), A(\. ) : \T \ra G ,
\text { where } C \in C^{1}(\T ,G)
\endeq
Then, $LA(\. ) : \T \ra g$ is a constant iff $A(\. ) $ is a torus morphism.
H. Eliasson's condition as he stated it is writen, in our notation,
\begin{equation} \label{condition UE}
\liminf \frac{| \langle LB_{i}, LB_{i+1}\rangle |}{|LB_{i}|.|LB_{i+1}|} <1
\end{equation}
i.e. the reductions of resonant modes of consecutive steps of the
algorithm form an angle bigger than some $\d >0 $, fixed throughout
the scheme, for an infinite number of steps. Since $LB_{i}$ commutes with
$A_{i}$, an equivalent statement is, by simple linear algebra, that
\begeq
\limsup \frac{d( A_{i+1} , \NN _{i} )}{d( A_{i+1}  , \Id )} > 0
\endeq
Here, $d ( \. ,\NN _{i} )$ is the Riemannian distance
from the normalizer of the torus of matrices simultaneously diagonalizable with $A_{i}$,
and $d$ is the Riemannian distance on $G$.

Eliasson's condition says, therefore, that a cocycle is UE provided that the matrix
$D_{i}$ of lemma \ref{lem est diag} is $\delta$-away from $\NN _{i}$, for an infinite
number of steps, for some $\d >0$.
%
%
This is equivalent to
\begin{equation} \label{eq rate approx e-F}
\limsup \big| \frac{\hat{F}_{i}(k_{i})}{\epsilon _{i}} \big| > 0
\end{equation}
We remark already that this condition implies that the corresponding small divisor of eq. \ref{local eq in Fourier}
is comparable with $|\hat{F}_{i}(k_{i})|$. Since the resonant modes $\hat{F}_{i}(k_{i})$ are (essentially) the Fourier
coefficients of a $C^{\infty }$ mapping, $|\hat{F}_{i}(k_{i})| = O(N_{n_{i}}^{-\infty})$, which implies that
under the condition of eq. \ref{condition UE},
\begeq
|\epsilon _{i}|_{\Z }= O(N_{n_{i}}^{-\infty})
\endeq
(at least along a
subsequence of $\{ n_{i} \}$), so that the estimation
$|\epsilon _{i}|_{\Z } \leq N_{n_{i}}^{-\nu }$ is very
optimistic in this case. In other words, a cocycle is UE provided that
the rotation around which we linearize at an infinite number of steps
of the K.A.M. scheme is (exactly resonant or) Liouville with respect to
$\a $. The genericity of condition of eq.
\ref{condition UE} is now clear.

Finally, we discuss the rareness of reductions of resonances, following \cite{El2002a} and
give the following lemma.
\begin{lem}
Along a subsequence of steps $n_{i_{j}}$ for which the condition of eq. \ref{condition UE} is
satisfied, resonant steps occur rarely:
\begeq
n_{i_{j}+1}-n_{i_{j}} \ra \infty
\endeq
\end{lem}
\begin{proof}
By construction,
\begeq
d(A_{i_{j}+1},\Id) = \sqrt{(\epsilon _{i_{j}})^{2}+|\hat{F}_{i_{j}}(k_{i_{j}})|^{2}}
\endeq
Since the condition of eq. \ref{eq rate approx e-F} is equivalent to eq. \ref{condition UE},
we obtain that, along the subsequence $i_{j}$,
\begeq
|\epsilon _{i_{j}} |<  \delta  |\hat{F}_{i_{j}}(k_{i_{j}})|
\endeq
for some fixed $\d >0$. Since $|\hat{F}_{i_{j}}(k_{i_{j}})| = O(N_{i_{j}}^{-\infty})$ and
since $\a \in DC$, we immediately obtain that
\begeq
k_{i_{j}+1}=O(N_{i_{j}}^{\infty})
\endeq
which ends the proof.
\end{proof}

If the condition for UE is not satisfied, we can only obtain the following estimate.
The constant $A_{i}$ is $K_{n_{i}}^{-1} = N_{n_{i}}^{-\nu}$-close to the identity, so that
the resonant mode $k_{i+1}$ must satisfy
\begin{equation*}
\frac{\gamma ^{-1}}{|k_{i+1}|^{\tau }} \leq N_{n_{i}}^{-\nu}
\end{equation*}
which implies $|k_{i+1}| \geq \gamma ^{-1/ \tau } N_{n_{i}}^{\nu / \tau } $. Since
$|k_{i+1}| \leq N_{n_{i+1}}= N_{n_{i}} ^{(1+\s )^{n_{i+1}-n_{i}-1}}$, we
obtain that
\begin{equation} \label{rareness of red of res}
n_{i+1} - n_{i} \gtrsim \frac{\log \nu  -\log \tau }{ \log (1+\s ) }
\end{equation}
where we remind that $\s $ is chosen positive and small. The conclusion of the lemma can
be obtained formally from this calculation by letting $\nu \ra \infty $.

\subsection{A concrete example}

The K.A.M. normal form of the cocycle also allows us to describe
quite explicitly a UE cocycle. We remind that in our notation,
$\RR (\pi /2)$ conjugates $\{ 1,0 \} _{g}$ to $\{ 0,1 \} _{g}$, and remark that we will be
using notation from continued fraction expansion (cf. \S \ref{Arithmetics, continued
fraction expansion}). Fix, now
$\{ q_{n_{i}} \}$ a sub-sequence of good denominators of $\a $, with
$n_{i}$ increasing fast enough. We omit the notation of the
sub-sequence. Then, the cocycle
\begin{equation*}
(\a , \{ e^{2i\pi q_{1} \a } ,0 \} _{G} . e^{\{ 0, \beta _{2} e^{2i\pi 2q_{1} \.} \} _{g}} )
\end{equation*}
is resonant for the first step of K.A.M. The reduction of the resonant mode by
$B_{1} = \{ e^{-2i\pi q_{1} \. } ,0 \}$ conjugates it to the resonant constant
\begin{equation*}
(\a , \RR (-\pi /2) \{ e^{2i\pi q_{2} \a} ,0 \} \RR^{*} (-\pi /2))
\end{equation*}
Therefore, if we perturb the initial cocycle by
\begin{equation*}
\exp (Ad(\RR (\pi /2). B_{1} (\. ) \RR ^{*}(\pi /2)).\{ 0, \beta _{3} e^{2i\pi 2q_{2} \.} \} _{g} )
\end{equation*}
it will be resonant in the next step. Then, iterating the construction we can obtain a cocycle satisfying the
smallness conditions of the K.A.M. scheme at all steps, as well as
Eliasson's condition for UE, since at all steps $\epsilon _{i} = 0$ (compare with eq.
\ref{eq rate approx e-F}).

In order to obtain a more general UE cocycle, one can replace $\RR (\pi /2)$ by a sequence
$\{ D_{i} \} $ of matrices such that $D_{i}$ is away from
all matrices simultaneously diagonalizable with $A_{i}$, and $q_{i}$
by a fast increasing sequence $\{ k_{i} \}$ such that
$|k_{i}\a |_{\Z }$ decreases exponentially fast. Such cocycles, as
well as typical cocycles close to them, are UE.

\subsection{Proof of theorem \ref{finite reduc thm}} \label{loc dens of fin dif}

In this section, we provide the proof of theorem \ref{finite reduc thm}.

We consider a cocycle $(\a , A_{1}e^{F_{1}(\. )}) \in \WW $
given in K.A.M. normal form and fix $0 \leq \s_{0} \leq \infty$.

The case where $\s_{0} = \infty$ corresponds to a smoothly reducible
cocycle, and has been proved already in \cite{KrikAst}. In our context,
the proof for this case reads as follows. By the almost reducibility
theorem, \ref{thm local a.r.} and lemma \ref{preliminary reduction}, the
sequence of conjugations $H_{n_{i}}$ constructed by the KAM scheme only
eliminates resonances and moreover satisfies the estimate of eq.
\ref{eq rate conv}. Consequently, if we write
\begeq
A_{1}e^{F_{1}(\. )} = H_{i}(\. +\a).A_{i}e^{F_{i}(\. )}.H_{i}^{*}(\.)
\endeq
then we immediately get
\begeq
A_{1}e^{F_{1}(\. )}.Ad(H_{i}(\.)).e^{-F_{i}(\. )} = H_{i}(\. +\a).A_{i}.H_{i}^{*}(\.)
\endeq
Since, by eq. \ref{eq rate conv}, the factor $Ad(H_{i}(\.)).e^{-F_{i}(\. )}$
can be made arbitrarily small by letting $i \ra \infty$, density of
smoothly reducible cocycles in $\WW$ is proved. We note that we showed that
the set of cocycles whose KAM normal form has a finite number of resonances
is dense in $\WW$. Such cocycles are clearly $\sm $ reducible.

We now concentrate on the case $\s_{0}<\infty$. A cocycle that is reducible
via a finite regularity transfer function has clearly an infinite normal
form. If the normal form of $(\a , A_{1}e^{F_{1}(\. )})$ is a finite product of
resonant modes, the cocycle can be perturbed so that the product becomes
infinite in the following way.
\begin{lem} \label{lem infinite KAM form}
Cocycles with infinite normal form are generic in $\WW$.
\end{lem}
\begin{proof}
If the cocycle has finite normal form, i.e. if for some $i_{l}$
\begeq
A_{1}e^{F_{1}(\. )} = H_{i_{l}}(\. +\a).A_{i_{l}}.H_{i_{l}}^{*}(\.)
\endeq
then one can perturb $A_{i_{l}}$ so that its eigenvalues become resonant.
One can then introduce a perturbation $F_{k_{l}}(\. )$ spectrally supported
on the resonant Fourier mode, and iterate the construction. Genericity is
straightforward.
\end{proof}
By Eliasson's theorem, \ref{thm El UE}, for a cocycle with an infinite
normal form to be reducible via a finite-regularity transfer function,
one has to impose some extra conditions.

This condition can be stated in terms of the sequence of angles $\theta$, as
shows the following lemma concerning the principal part of a cocycle in
KAM normal form after $i$ steps of the KAM scheme.
\begin{lem} \label{conj in fin reg iter}
Consider a resonant cocycle, perturbed only by the resonant Fourier mode:
\begin{equation*}
(\a , A_{i}e^{\tilde{F}_{i}(\. )}) =
(\a ,
\{ e^{2i\pi (k_{i}\a + \epsilon _{i})} , 0 \} .
\exp\{ 0 ,\hat{F}_{i} (k_{i}) e^{2i\pi k_{i})\.}  \})
\end{equation*}
Such a cocycle can be reduced by a conjugation $G'_{i}(\.)$ satisfying
the estimate
\begin{equation*}
\|G'_{i}(\.)\|_{s} \simeq |k_{i}|^{s} |\sin \theta _{i}/2|
\end{equation*}
where
\begeq
\theta _{i} = \arctan \frac{|\hat{F}_{i}(k_{i})|}{\epsilon _{i}}
\endeq
\end{lem}
The estimate for the conjugation is optimal. Of course, reducing conjugations exist
satisfying worse estimates.
\begin{proof}
We recall that
\begin{equation*}
Conj_{B_{i}(\. )}(\a , A_{i}e^{\tilde{F}_{i}(\.)}) = (\a , A_{i+1})
\end{equation*}
If we now let $D_{i}$ diagonalize $A_{i+1}$ in the basis where
$A_{i}$ is diagonal as in lemma \ref{lem est diag}, then $D_{i}$ is
$\sin \theta _{i}/2$-away from the identity in $SU(2)/\TT _{i}$.

Then, the conjugation
$G^{\prime}_{i} = B_{i}^{*}(\. )D_{i} B_{i}(\. )$ is $N_{i}^{\s_{0}}|\sin \theta _{i}|$-away
from the $\Id$ in the $H^{\s_{0}}$ norm, since it has the
explicit form
\begeq
G_{i}'(\. ) = 
\{\cos \theta _{i} /2 , i e^{2i\pi (k_{i}\. + \f_{i})} \sin \theta _{i} /2 \}_{SU(2)}
\endeq
This conjugation thus satisfies the estimate announced in the lemma and acts as follows:
\begeq
\begin{array}{r@{}l}
Conj_{G'_{i}(\. )}(\a , A_{i}e^{\tilde{F}_{i}(\.)}) &= Conj_{B_{i}^{*}(\. )D_{i}}(\a , A_{i+1})
\\
&= Conj_{B_{i}^{*}(\. )D_{i}}(\a , \{ e^{2i\pi \epsilon _{i}} , 0 \}_{G} . e^{\{ 0 , \hat{F}_{i}(k_{i}) \}_{g}})\\
&= Conj_{B_{i}^{*}(\. )}(\a , \{ e^{2i\pi \r _{i}} , 0 \}_{G})\\
&= (\a , \{ e^{2i\pi (\r _{i}+k_{i}\. \a )} , 0 \}_{G})
\end{array}
\endeq
where
\begin{equation} \label{eq def rho}
\r _{i} = \sqrt{(2\pi \epsilon _{i})^{2} + |\hat{F}_{i}(k_{i})|^{2}}
\end{equation}

Optimality of the estimate is an immediate consequence of the corresponding linear
problem, following the proof of lemma \ref{loc conj lem}, and especially of the
estimation of the small denominator in eq. \ref{local eq in Fourier}.
\end{proof}

We can now prove theorem \ref{finite reduc thm}, using the notation introduced in lemma
\ref{conj in fin reg iter}.
\begin{proof} [Proof of theorem \ref{finite reduc thm}]
Fix $\s _{0} \in [0,\infty )$.
Direct calculation shows that, up to the first order,
\begin{equation*}
Conj_{G_{i}^{\prime}}(\a , A_{i}e^{F_{i}}) = (\a ,
\{ e^{2i\pi (k_{i+1}\a + \epsilon _{i+1})} , 0 \} .
\exp\{ 0 ,\hat{F}_{i+1} (k_{i+1}) e^{2i\pi (k_{i+1}+k_{i})\.}  \})
\end{equation*}
We can therefore continue the reduction by fixing $i_{0}$, a resonant step, and then
defining
\begeq
\begin{array}{r@{}l}
\tilde{B} _{i}(\. ) &= \exp \{ 2i\pi \sum _{j=i_{0}}^{i} k_{j} \. , 0 \} \\
D_{i} &\in  SU(2) / \TT _{i_{0}}  \text{ such that}\\
D_{i} A_{i} D_{i}^{*} &\in  \TT _{i_{0}} \\
G_{i}^{\prime} (\. ) &= \tilde{B} _{i}^{*}(\. ) D_{i} \tilde{B} _{i}(\. ) \\
H_{i}^{\prime} (\. ) &=  \prod _{j=i_{0}}^{i} G_{j}^{\prime} (\. )
\end{array}
\endeq
where $\TT _{i_{0}}$ is the torus of matrices simultaneously diagonalizable
with $A_{i_{0}}$. Then, by the estimate in lemma \ref{conj in fin reg iter},
$H_{i}^{\prime} (\. ) \ra H(\. )$ in $H^{\s_{0}}$ iff
\begeq
\sum _{i\geq i_{0}}|k_{i}|^{s} |\sin \theta _{i}/2| < \infty
\endeq
Since $k_{i}$ is of the same order as $N_{i}$ and $\sin \w \simeq \w $ for $\w$
small, this condition on the sequence $\theta$ is equivalent to
\begin{equation} \label{cvgence cdtion}
\{ \theta _{i }\} = \left\lbrace \arctan \frac{\| F_{i}\| _{L^{2}}}{| \tilde{\epsilon} _{i} |}
\right\rbrace _{i \geq i_{0}} \in \tilde{h}^{\s_{0}}
\end{equation}
where the weights defining the $h^{s}$ norm are $N_{i}^{s}$ (cf. eq. \ref{def h spaces}).

Summarizing, we have obtained that a cocycle is $H^{\s_{0}}$ reducible and not reducible
in any higher regularity iff the condition of eq. \ref{cvgence cdtion} is satisfied.

For fixed $\s _{0}$, this condition is an $F_{\s}$ condition with empty interior in
$\WW$ as shows the following argument.

Firstly, the cocycles having an infinite KAM normal form form a dense, actually a $G_{\d }$,
subset of $\WW$.

Secondly, any cocycle can be perturbed in order to satisfy such a
condition, by letting $i_{0}$ be large enough, truncating the KAM normal form and appending
a normal form starting at the index $i_{0}$ and satisfying the condition of eq.
\ref{cvgence cdtion}. An example of such a KAM normal form can be constructed iteratively
by hand (both points by lemma \ref{lem infinite KAM form}).

Finally, the condition
\begeq
\left\lbrace \frac{\| F_{i}\| _{L^{2}}}{| \epsilon _{i} |}
\right\rbrace _{i \geq i_{0}} \in h^{\s_{0}}
\endeq
where we do not impose optimality of $\s_{0}$ is an $F_{\s}$ condition in the space of
real sequences. The condition of eq. \ref{cvgence cdtion} where optimality is imposed is
an additionnal $F_{\s}$ condition.
\end{proof}

A simple perturbation argument allows us to obtain the following
corollary (cf. \cite{Chav12}).
\begin{cor}
Let $\a \in DC (\gamma , \tau )$, and $\s _{0} \geq s _{0}$, where
$s_{0}$ is given by the smallness condition of the local almost
reducibility theorem. Consider the set of $C^{\s _{0}}$-smooth cocycles
over the diophantine rotation $\a \in \T ^{d}$, denoted by
$SW^{\s_{0}}_{\a } (\T ^{d} , SU(2))$ and $\WW$, the neighborhood of
constants as in theorem \ref{finite reduc thm}. Then,
$C^{\s_{0}-d}$-reducible cocycles are $C^{\s_{0}-d}$ dense in $\WW$.
\end{cor}
\begin{proof}
We perturb the given cocycle in $C^{\s_{0}}$ so that it becomes
$C^{\infty}$, and we apply theorem \ref{finite reduc thm}.
The Sobolev injection theorem accounts for the loss of derivatives.
\end{proof}

\section{DUE in $SW_{\a}^{\infty}(\T ,SU(2)) $} \label{section proof main thm}


\subsection{The main ideas}

In this section, we prove theorem \ref{thm DUE}. The strategy of the proof is to apply
proposition \ref{prop DUE qp}. The approximants $f_{n}$ of the proposition are provided
by the almost reducibility theorem, thm. \ref{thm local a.r.}. We will then prove that
under a condition of the flavour of the one of eq. \ref{condition UE}, proposition
\ref{prop DUE qp} can be applied, thus proving that a generic cocycle in the KAM regime
is DUE.

We think that, in view of the normal form obtained in lemma \ref{preliminary reduction}, the proof of local
genericity of DUE also clarifies to a certain extent the proof of local genericity of unique ergodicity.

In a nutshell, the proof consists in
using the K.A.M. scheme of the local almost reducibility theorem in order to solve
\begin{equation*}
\psi \circ \Phi - \psi = \phi - \int \phi
\end{equation*}
up to an arbitrarily small error in $C^{\infty} (\T \times SU(2))$, where $\phi \in C^{\infty} $
is a given function (from now on of $0$ mean value) and
$\Phi = \Phi _{1} = (\a ,A_{1}e^{F_{1} (\. )}) \in SW^{\infty}_{\a}$ is a given cocycle close to constants
satisfying a generic condition. In order to simplify the proof, we use the following
observation.
\begin{observation}
Let $\pi _{m}$ be an irreducible representation of $SU(2)$. Then, the space
$L^{2}(\T ^{d}) \times \EE _{\pi _{m}}$ is invariant under the action of
$SW(\T^{d} ,SU(2))$.
\end{observation}
In particular, since fibered conjugations act as cocycles with $0$ as their frequency,
we obtain the following.
\begin{observation}
The space
$L^{2}(\T ^{d}) \times \EE _{\pi _{m}}$ is invariant under the action of
conjugations.
\end{observation}
We can therefore assume that $\f $ is spectrally supported in a single irreducible
representation of $SU(2)$, .i.e. that
all harmonics belong to a fixed space $\EE _{\pi _{m}}$, $m>0$:
\begin{equation*}
\f (x, z,w) = \sum _{\substack{k \in \Z ^{*}\\
0 \leq j,p \leq m }}
\f ^{m}_{k,j,p} e^{2i\pi k x}\p _{m}^{j,p}(z,w )
\end{equation*}
We remind our assumption that that the cocycle $\Phi$ is given in its KAM normal form.

The cohomological equation can be solved by linearizing with respect to the dynamics
(thus introducing an error term), and solving for a finite truncation
(both in frequencies and the harmonics).
The linearized equation has a large number of obstrutions, the functions
in the kernels of distributions preserved by a constant (moreover,
resonant) cocycle. We solve, therefore,
\begin{equation*}
\psi _{1} \circ \Phi _{1} - \psi _{1} - \f = R_{1} \f + Ob_{1} \f +E_{1} = \f _{2}
\end{equation*}
where $R_{1}$ is a rest operator, $Ob_{1}$ the projection on the obstructions, and 
$E_{1}=O_{C^{s}}(\| \partial \psi _{1} \|_{s} \|F _{1} \|_{s} )$, for every $s$.
When we make one step of the K.A.M. scheme and introduce the resulting change of coordinates we obtain the cocycle
$\Phi _{2} =(\a ,A_{2}e^{F_{2} (\. )}) $ which is a perturbation of a different
constant, and the order of the perturbation
is much smaller. We now obtain the equation
\begin{equation*}
\psi _{2}^{\prime} \circ \Phi _{2} - \psi _{2}^{\prime} = \f _{2}^{\prime}
= R_{1} ^{\prime} \f + Ob_{1}^{\prime} \f + E_{1}^{\prime}
\end{equation*}
where the prime indicates that quantities are expressed in the new coordinates. Since
$F_{2}$ is quadratic with respect to $F_{1}$, the linearization error
term will be smaller, allowing more frequencies and harmonics in the truncation (and a larger space of obstructions).
As for $ Ob_{1}^{\prime} \f $, they are constant with respect to $x$.
If, now, $A_{2}$ does not commute with $A_{1}$,
lemma \ref{change of coord harm} says that only a fraction (in the $L^{2}$ norm) of these
obstructions persists as an obstruction to the new equation. We thus solve
\begin{equation*}
\psi _{2}^{\prime} \circ \Phi _{2} - \psi _{2}^{\prime} = Ob_{2}^{\prime} \f _{2}^{\prime} +
R_{2} \f_{2} ^{\prime} + E_{2}^{\prime}
\end{equation*}
When we inverse the change of coordinates, we find that
$ Ob_{2} Ob_{1}^{\prime} \f$, in the standard coordinates,
has its spectral support in the same ball as $Ob_{1}$, but a
smaller $L^{2}$, and therefore $H^{s}$, norm.
When we change coordinates following the $3$rd step of K.A.M., $ (Ob_{2} Ob_{1}^{\prime})^{(2)} \f$
will still not depend on $x$ in the new coordinates, and $ Ob_{3}(Ob_{2} Ob_{1}^{\prime})^{(2)} \f$ will have a smaller
$L^{2}$ norm for the same reason. Inversion of the change of coordinates, however, will create new frequencies.
We need, therefore, to impose a condition guaranteeing that the spaces of obstructions of successive constants are almost orthogonal infinitely often
for every fixed $m$. This is the refinement of condition \ref{condition UE} needed in
order to ensure DUE.

We thus obtain that the closure of coboundaries coincides with $0$ average functions, which is equivalent to
saying that the only distribution preserved by the dynamics is the product of the Haar measures.

The proof admits the following dual interpretation. The K.A.M. scheme, when the product of
conjugacies diverges, can be seen as a formal solution to the equation (with respect to $H_{\infty}(\. )$)
\begin{equation*}
Conj_{H_{\infty}(\. )}(\a ,A e^{F(\. )}) = (\a , \Id )
\end{equation*}
This is so, because the sequence of constants constructed by the scheme then converges to the $\Id$ and
the perturbation converges
to $0 $ in $C^{\infty}$. The condition for DUE assures that, when we invert the scheme
in order to obtain the original cocycle from its formal limit, the Fourier coefficients
of all the invariant distributions are pushed to infinity by the conjugations reducing
the resonant modes.

Cocycles preserving only a finite number of distributions in the border
between DUE and UE do not exist,
since if we let the condition implying DUE degenerate, distributions
of arbitrarily high order are created, and no linear dependence
between their Fourier coefficients can be forced. Consequently, the space of coboundaries
passes from a codimension $1$ to an infinite-codimension space in $\sm (\T \times G, \C )$.

\subsection{The basic lemmas}
We suppose that the cocycle is given in normal form. Since the K.A.M. scheme produces only conjugations far from the
identity, we study the behavior of harmonics on $\T \times G $ under such conjugations. The
proof is by direct calculation.
\begin{lem} \label{transformation of harmonics}
If $B(\. ) : \T \ra SU(2) $ is of the form
\begin{equation*}
B(\. )=
\begin{pmatrix}
e^{2i\pi k_{r} x} & 0 \\
0 & e^{-2i\pi k_{r} x}
\end{pmatrix}
\end{equation*}
then the harmonic $e^{2i\pi k x} \p _{m}^{j,p} (z,\bar{z},w,\bar{w})$ is mapped to
\begeq
e^{2i\pi (k+(m-2p)k_{r}) x}\p _{m}^{j,p}(z,\bar{z},w,\bar{w})
\endeq
\end{lem}
The definition of the harmonic functions $\pi $ can be found in eq. \ref{eq def pi}.
We remark that in the proof of almost reduciblilty, cf. lemma \ref{loc conj lem}, we saw that
$2$-periodic conjugations might be applied. They are of the form
\begin{equation*}
\begin{pmatrix}
e^{2i\pi k_{r} x/2} & 0 \\
0 & e^{-2i\pi k_{r} x/2}
\end{pmatrix}
\end{equation*}
and they are actually an artefact of the method, as we proved in \cite{NKPhD}, \S $9$.
In the context of the present article, in order to account for the use of such conjugations
we should make the cocycle act on fuctions on $2\T \times SU(2)$, where $2\T = \R / 2\Z $.
This would only complicate notation, and otherwise lead to the same conditions, as inspection
of the rest of the proof can show. We therefore ignore this technicality in order to keep
notation simpler.

As a consequence of the previous lemma, we obtain the following important estimate, implying that,
after applying such a conjugation, the $C^{s}$ norms of harmonics deteriorate by factors
of the order of $k_{r}^{ s}$.
\begin{lem}
Let $\f (\. ) \in C^{\infty} (\T \times SU(2)) $, $B(\. ) $
as in the previous lemma, and call
\begin{equation*}
\tilde{\f } (x, \{ z,w \} )  = \f (x, B^{-1}(x ).\{ z,w \} )
\end{equation*}
Then
\begin{equation*}
\| \tilde{\f } \| _{s} \leq C_{s} (k_{r})^{s+2} \| \f \| _{0} 
\end{equation*}
\end{lem}

\begin{proof}
The proof follows the lines of similar estimations, using the expansion
$\f (x,\{ z,w \} )= \sum _{k \in \Z , 0 \leq j,p \leq m } \f ^{m}_{k,j,p} e^{2i\pi kx} \p _{m}^{j,p} $
\begin{eqnarray*}
\| \tilde{\f } \| _{s}^{2} & \leq & C_{s} \sum _{k,j,p}
(1+k^{2} + m^{2} )^{s + 2}  | \tilde{\f } ^{m}_{k,j,p} |^{2} \\
& = & C_{s} \sum (1+(k + (m-2p)k_{r} )^{2} + m^{2} )^{s + 2} | \f ^{m}_{k,j,p} |^{2} \\
& \leq & C_{s} k_{r} ^{2(s+2)} \sum (1+k^{2} + m^{2} )^{s + 2} | \f ^{m}_{k,j,p} |^{2} \\
& \leq & C_{s} k_{r} ^{2(s+2)}  \| \f \| _{0}^{2}
\end{eqnarray*}
where as usual the constant may change from one line to another.
\end{proof}

Let us now return to the solution of the cohomological equation.
Since our strategy is to use almost reducibility through proposition \ref{prop DUE qp},
the following elementary lemmas will be useful.

\begin{lem}
Let $(\a ,A)$ be a constant cocycle. Then, each function of the form
$e^{2i\pi k x}\p _{m}^{j,p}(z,w )$ is an eigenfunction under composition by the cocycle.
\end{lem}
\begin{proof}
By direct calculation, we obtain
\begeq
\begin{array}{r@{}l}
e^{2i\pi k x}\pi _{m}^{j,p}(z,w ) \circ (\a ,A) &= e^{2i\pi k (x+\a)}\p _{m}^{j,p}(e^{2i\pi a}z,w )\\
&=e^{2i\pi (k\a +(m-2p) a)}e^{2i\pi k x}\p _{m}^{j,p}(z,w )
\end{array}
\endeq
where $A=\{ e^{2i\pi a},0 \}_{G}$.
\end{proof}
We immediately obtain the following lemma.
\begin{lem}  \label{lem cohom obs}
Let $(\a ,A)$ be a resonant cocycle,  i.e. suppose that in the notation of the previous
lemma $a = k_{r}\a$, for some $k_{r } \in \Z ^{*}$. Then, a function $\f $ of the form
$e^{2i\pi k x}\p _{m}^{j,p}(z,w )$ is not a coboundary over $(\a ,A)$ iff the indices
$k,j,p$ satisfy the resonance relation
\begin{equation} \label{cohom resonances}
k + k_{r}(m-2p) = 0
\end{equation}
\end{lem}
\begin{proof}
The cohomological equation over a constant cocycle $(\a , A)$ reads
\begin{equation*}
\psi (x +\a , A^{*}. \{ z,w \} ) - \psi (x , \{ z,w \} ) =\f  (x , \{ z,w \} ) 
\end{equation*}
In terms of coefficients, and for $A = \{e^{2i\pi a}, 0 \}$, we have
\begin{equation*}
(e^{2i\pi (k\a +(m-2p) a)}-1)\psi  ^{m}_{k,j,p} = \f  ^{m}_{k,j,p}
\end{equation*}
Thus, the resonant modes (the obstructions of the equation), are
given by the relation $k\a +(m-2p) a \in \Z$.
\end{proof}
The Diophantine condition imposed on $\a$, cf. \ref{def DC}, implies the following lemma.
\begin{lem} \label{lem cohom res 1 mode}
Let $(\a ,A)$ be a resonant cocycle. Then, for a function $\f  $ of the form
$e^{2i\pi k x}\p _{m}^{j,p}(z,w )$ whose indices $k,j,p$ are non-resonant, i.e.
\begeq
k + k_{r}(m-2p) \neq 0
\endeq
there exists a solution of the cohomological equation, $\psi$, satisfying the estimate
\begeq
\|\psi \| _{s} \leq C_{s} |(k + k_{r}(m-2p))\. \a|_{\Z}^{\t + 2}\|\f  \| _{s}
\endeq
\end{lem}

In view of the K.A.M.
normal form, the resonance condition is implied by
\begin{equation*}
a = k_{r}\a \in \a \Z \text{  and } k +(m-2p)k_{r} =0
\end{equation*}
This relation is, of course, valid in a coordinate system $\{ z,w \} \in G $ where $A$ is diagonal.
For a fixed vector $e^{2i\pi k x} z^{p} w^{m-p}$ (where $m,p$ are fixed and $k \in \Z$) there is at most one
solution to the equation \ref{cohom resonances}.

The cohomological equation for the perturbed cocycle reads
\begin{equation*}
\psi (x +\a , e^{-F(\. )}A^{*}. \{ z,w \} ) - \psi (x , \{ z,w \} ) =\f  (x , \{ z,w \} ) 
\end{equation*}
which we rewrite as
\begin{equation*}
\psi  (x +\a , A^{*}. \{ z,w \} ) - \psi  (x , \{ z,w \} ) =\f   (x , \{ z,w \} )+
O(\| \partial \psi \|_{s} \| F \|_{s})
\end{equation*}
where the last term is a linearization error term. This fact obviously poses a constraint on the number of
modes that can be reduced, relatively to the size of the perturbation.
The resonance relation holds approximately in this setting,
since the cocycle is a perturbation of a resonant constant, $(\a , A e^{F(\. )} )$ where
$a = k_{r}\a + \epsilon $,
with $|\epsilon | < |k_{r}|^{-\nu } $. We have $ |\hat{F} (k_{r})| \ll |k_{r}|^{-1} $,
and our condition (condition \ref{condition DUE}) will imply that $|\epsilon | \lesssim |\hat{F} (k_{r})|$
infinitely often in the K.A.M. scheme. Therefore, the resonant modes,
for which the small denominator is $\ll |k_{r}|^{-1}$, are
\begin{equation*}
k +(m-2p) k_{r} =0 , 0\leq p \leq m
\end{equation*}

We recall the following lemma from \cite{El2002a} and \cite{NKPhD}.
\begin{lem} \label{Dioph csts}
Let $\a \in DC(\gamma ,\tau )$, $a \in \T $, $N \in \N ^{\ast }$ big enough,
and $K = N ^{\nu } >0$, with $\nu > \tau $, so that $K\geq 2^{\tau +1}\gamma N^{\tau }$. Then there exists at most one
$k=k_{r} $, $0\leq |k_{r}| \leq N $, such that
\begin{equation*}
\left\vert a -k \a \right\vert _{\Z }<K^{-1},~0<|k|\leq N
\end{equation*}
Moreover, if such a $k_{r}$ exists and $\left\vert a -k_{r} \a \right\vert _{\Z } = \epsilon$
\begin{equation*}
|a -k\a  |_{\Z} \geq  K^{-1} \text{ for } 0<|k- k_{r} |\leq N^{\prime }
\end{equation*}
where $N^{\prime } = \left( \frac{\gamma}{1+K\epsilon} \right) ^{1/ \tau } N^{\nu / \tau} \geq 
\left( \frac{\gamma}{2} \right) ^{1/ \tau } N^{\nu / \tau} $.
\end{lem}
Therefore, if the resonance relation is satisfied with an $\epsilon $ of room, as in the lemma, and if we call
$T^{k_{r}}_{N '}$ the operator
acting on $C^{\infty} (\T \times G)$ and projecting on the modes $(k,j,p,m)$ such that
$0<|k- (m-2p) k_{r} |\leq N^{\prime }$,
\begin{equation*}
T^{k_{r}}_{N'} \f   (x , \{ z,w \} ) = \sum _{0<|k+ (m-2p)k_{r} |\leq N'}
\f  _{k,j,p}^{m} e^{2i\pi k x} \p ^{j,p}_{m}
\end{equation*}
the denominators in the image of $T^{k_{r}}_{N'}$ are $\geq N^{-\nu } $. On the
other hand, if $Ob ^{k_{r}}_{N'} $ is the projection on the resonant modes,
\begin{equation*}
Ob ^{k_{r}}_{N'} \f   (x , \{ z,w \} ) =
\sum _{
\substack {k+ (m-2p)k_{r}=0 \\
|k |\leq N' }
}
\f  _{k,j,p}^{m} e^{2i\pi k x} \p ^{j,p}_{m}(z,w)
\end{equation*}
the denominators in the image of $Ob^{k_{r}}_{N'} $ are $< N^{-\nu }$,
and the corresponding modes cannot be eliminated with good
estimates. We thus obtain the following lemma.
\begin{lem}
Let $(\a ,A)$ be a resonant cocycle. Then, for a function $\f  $ spectrally supported in the
image of $T^{k_{r}}_{N'}$, i.e. such that
\begeq
T^{k_{r}}_{N'} \f  = \f 
\endeq
there exists a solution of the cohomological equation, $\psi$, satisfying the estimate
\begeq
\|\psi \| _{s} \leq C_{s} N^{s+\n + 2}\|\f   \| _{0}
\endeq
\end{lem}
\begin{proof}
The proof follows directly by combining lemmas \ref{lem cohom res 1 mode} and
\ref{Dioph csts}.
\end{proof}

Therefore, the equation
\begin{equation*}
\psi  (x +\a , A^{*}. \{ z,w \} ) - \psi  (x , \{ z,w \} ) = T^{k_{r}}_{N'} \f   (x , \{ z,w \} )
\end{equation*}
admits a solution that has the same harmonics as $T^{k_{r}}_{N'} \f   (x , \{ z,w \} )$.
Consequently, if the rest operator $R^{k_{r}}_{N'}$ is defined so that
\begin{equation} \label{eq def rest op}
T^{k_{r}}_{N'} + Ob^{k_{r}}_{N'} + R^{k_{r}}_{N'} = \Id
\end{equation}
we have obtained a solution to the equation
\begin{equation*}
\psi  \circ (\a , A) - \psi = \f  - Ob^{k_{r}}_{N'} \f  -
R^{k_{r}}_{N'} \f 
\end{equation*}
where of course, the obstruction is of the order of $\f  $, and the rest is small:
\begin{eqnarray*}
\| Ob^{k_{r}}_{N'} \f  \|_{s} &\leq & C_{s} (N')^{s+2} \| \f  \| _{0} \\
\| R^{k_{r}}_{N'} \f  \|_{s} &\leq & C_{s,s'} (N')^{s-s'+2} \| \f  \| _{s'}
\end{eqnarray*}

Since the cohomological equation is linear, we can obtain the following lemma,
connecting the solution we have constructed with the KAM normal form.
\begin{lem} \label{lem cohom eq w pert}
If $\Phi = (\a ,A^{F(\. )})$ is in normal form, $A \in G$ is as in lemmas \ref{lem cohom res 1 mode} and
\ref{Dioph csts},
and $\f  \in C^{\infty}_{0} $, then there
exists a function $\psi $ satisfying
\begin{equation*}
\psi \circ \Phi  - \psi = \f - Ob^{k_{r}}_{N'} \f    -
R^{k_{r}}_{N'} \f  - E
\end{equation*}
with the estimates
\begin{eqnarray*}
\| E \|_{s} &\leq &  C_{s} K (N') ^{s+4} \| \f  \| _{0} \| F  \|_{s} \\
\| R^{k_{r}}_{N'} \f  \|_{s} &\leq & C_{s,s'} (N')^{s-s'+3} \| \f  \| _{s'} \\
\| \psi \| _{s}  &\leq & N ^{s+\nu + 3} \| \f  \| _{0}\\
\| Ob^{k_{r}}_{N'} \f  \|_{s} &\leq & C_{s} (N')^{s+2} \| \f  \| _{0} 
\end{eqnarray*}
\end{lem}

\begin{proof}
Let
\begin{equation*}
\sum _{\substack{k \in \Z \\
 0 \leq j,p \leq m}}
\f  ^{m}_{k,j,p} e^{2i\pi k x}\p _{m}^{j,p}(z,w )
\end{equation*}
and apply lemma \ref{lem cohom res 1 mode} to
\begin{equation*}
\sum _{\substack{|k| \leq N' \\
0 \leq j,p \leq m }}
\f  ^{m}_{k,j,p} e^{2i\pi k x}\p _{m}^{j,p}(z,w )
\end{equation*}
in order to obtain $\psi$ satisfying
\begin{equation*}
\psi  \circ (\a , A) - \psi = \sum _{\substack{|k| \leq N  \\
0 \leq j,p \leq m }}
\f  ^{m}_{k,j,p} e^{2i\pi k \.}\p _{m}^{j,p} - Ob^{k_{r}}_{N'} \f 
\end{equation*}
as well as the estimate announced in the lemma. The obstruction term
$Ob^{k_{r}}_{N'} \f $, as well as the rest term $R^{k_{r}}_{N'}$ satisfy
the respective estimates.

Finally, the error term $E$ coming from the linearization can be estimated by
\begin{eqnarray*}
\| \psi (x +\a , A^{*}. \{ z,w \} ) - \f   (x+\a ,e^{-F(\. )}A^{*}. \{ z,w \} )  \|_{s}
&\leq & C_{s} (N')^{2} \| \psi \| _{s+1}  \| F  \|_{s} \\
&\leq & C_{s} K (N') ^{s+4} \| \f  \| _{0} \| F  \|_{s}
\end{eqnarray*}
This concludes the proof of the lemma.
\end{proof}

\subsection{Iterative construction of the approximate solution}
\label{sec constr approx sol}

We now relate the solution constructed in the previous section with the
KAM scheme and exploit almost reducibility.
It is in this point that a condition concerning successive constants, like condition of eq. \eqref{condition UE},
becomes important.

The reason is clearly
the fact that given a size of a perturbation, we can solve the cohomological equation only up to a certain
error, and by excluding the resonant harmonics which we have called obstructions. The error converges
to $0$ faster than the rate of growth of the conjugations, so the non-resonant part is easy to cope with,
thanks to the fast convergence of the K.A.M. scheme. In order to cope with the obstructions,
we need to impose a condition assuring that obstructions do not accumulate. They accumulate, for example, if
the cocycle is reducible, and the scheme will show an accumulation of obstructions if, for example, the
cocycle is $C^{\infty }$ reducible to a Liouvillean constant, in which case the scheme generically produces an infinite
number of reductions of resonances, due to the arithmetics of the dynamics in the fibers. We will see that
obstructions will not accumulate if the angles between sucessive reductions of resonances
visit the roots
of Legendre polynomials (cf. lemma \ref{Legendre roots}) infinitely often. This is stronger than imposing that
the constants around which we linearize in successive steps of K.A.M. do not commute as
in eq. \eqref{condition UE}.

The last lemma of the previous section, lem. \ref{lem cohom eq w pert}, takes care of the case
of a function supported in the low frequences and modes and having no obstructions. In that case,
we can obtain a first good approximation to writing the function as a coboundary. However, the set of
functions spectrally supported in the obstructions form an open dense set.

The next lemma takes care of that case. It states that, if a certain condition on the angles $\theta$
is satisfied, then some obstructions can be solved for after two steps of the KAM scheme. The
hypotheses on the cocycle just put the properties of a cocycle in normal form in the
context. We will use the notation introduced in the proof of lemma \ref{loc conj lem}.

\begin{lem} \label{iterative lemma}
Let $(\a , A_{1} e^{\tilde{F}_{1} (\. ) }e^{F_{1}^{\prime} (\. ) })$ be a cocycle and $\f  $ a function satisfying
\begin{enumerate}
\item \label{rec prop start} $A_{1} = \{e^{2i\pi ( k_{1}\a + \epsilon _{1} ) }, 0 \} \in G $ is a resonant constant, up to
$|\epsilon _{1}| \ll N_{1}^{-\nu }$, and the resonant mode
is $k_{1}$, with $|k_{1}| \leq N_{1}$ but of the order of $N_{1}$.
\item $\tilde{F}_{1} (\. ) = \{  0 , \tilde{F}_{1,z}(\. )\}_{g} : \T \ra g $ is small and $\s (\tilde{F}_{1,z}) = \{ 2k_{1} \}$.
\item $\f  : \T \times G \ra \C $ is spectrally supported in an obstruction of a given bidegree $(m-p,p)$,
and $m-2p = b$:
\begeq
\f (x,\{ z,w \} ) = Ob_{1}\f (x,\{ z,w \} ) = \sum _{\substack{
k+b k_{r} =0 \\
0\leq j \leq m
}} \f  ^{m}_{k,j,p} e^{2i\pi k x} \p _{m}^{j,p}(z,w )
\endeq
with $m \ll N_{1}$
\item $F_{1}^{\prime} (\. )$ is quadratic with respect to $\tilde{F}_{1} (\. )$ in the sense of the KAM scheme
\item \label{rec prop end} $\s (Ad(B_{1}(\. )).F_{1}^{\prime}) $ is outside some ball of
radius $\approx N_{2} \gg N_{1}$, where $B_{1}(\. )$ is the
conjugation reducing $(\a , A_{1} e^{\tilde{F}_{1} (\. ) }) $ to $(\a , A_{2} ) $ with
\begeq
A_{2} = \{e^{2i\pi \epsilon _{1}} ,0\}  .\{ 0 ,e^{\{ 0 , \hat{F}_{1} (k_{1}) \} }\}
\endeq
\item The cocycle
\begin{equation*}
Conj_{B_{1}(\. )}  (\a , A_{1} e^{\tilde{F}_{1} (\. ) }e^{F_{1}^{\prime} (\. ) }) =
(\a , A_{2} e^{\tilde{F}_{2} (\. ) }e^{F_{2}^{\prime} (\. ) })
\end{equation*}
satisfies conditions in items \ref{rec prop start} to \ref{rec prop
end}, with $N_{2} , N_{3}$ replacing $N_{1}, N_{2}$.
\item The angle between $Ob _{2} $ and $ Ob_{3} $ in $(\PP _{m})^{*} $ is close to a root of the corresponding Legendre
polynomial $p_{m/2}$ (cf. lemma \ref{Legendre roots}):
\begin{equation} \label{Legendre lemma}
\|Ob_{3} \circ Ob_{2} \| _{\mathcal{L} (L^{2})} \leq C  N_{2} ^{-s_{0}}
\end{equation}
for some constants $C>0$ and $s_{0} \gg 0$.
\end{enumerate}
We define the obstruction, truncation and rest operators for $(\a , A _{i}) $ as in the previous section. Then,
\begin{enumerate}
\item $\widetilde{Ob} _{1} \f  ( z,w ) =  Ob_{1} \f  \circ (0,B_{1} (\. ))  ( z,w )=
\sum _{k+bk_{1} = 0} \f  ^{m}_{k,j,p} \p _{m}^{j,p}  ( z,w ) $ does not depend on $x \in \T$
\item We have
\begin{equation*}
\|( Ob_{2} ( \widetilde{Ob} _{1} \f  )) ) \| _{L^{2}} \leq  \| \f  \| _{L^{2}}
\end{equation*}
\item $( Ob_{2} ( \widetilde{Ob} _{1} \f  )) ) \circ (0,B^{*}_{1}(\. ))$ is spectrally supported in the image of $Ob_{1}$
\item If we call $\widetilde{ Ob}_{2}  \f    =  ( Ob_{2} ( \widetilde{Ob} _{1} \f  )) ) \circ (0,B_{2}(\. ))$,
then $\widetilde{Ob}_{3} = Ob_{3} ( \widetilde{ Ob}_{2}  \f  ) ) \circ (0,B^{*}_{1}.B^{*}_{2}(\. ))) $ satisfies
\begin{equation*}
\|\widetilde{Ob}_{3} \| _{s}  \leq C_{s} N_{2}^{s-s_{0}+\l } \| \f  \| _{L^{2}} , \, 0 \leq s \leq s_{0}- \l
\end{equation*}
where $C$ is the constant in eq. \ref{Legendre lemma}, $C_{s} >0$ depends only on $s$, and $\l >0$ is a universal
constant.
\end{enumerate}
\end{lem}
Using this lemma, we can eliminate the obstructions activated at a certain step of the
K.A.M. scheme, up to a small error,
if we wait long enough so that the projections corresponding to two successive constants form an angle close to
a root of $p_{m/2}$. The lemma indicates the following condition under which DUE should
hold, see lemma \ref{Legendre roots} for the notation.
\begin{cond} \label{condition DUE}
The condition for DUE is that for any $m\in \N ^{*} $,
along a subsequence of $\{ n _{i} \}$,
\begin{equation*}
\lvert \cos ^{2} (\theta _{i} ) 
- \xi _{m/2} \rvert = O(N_{n_{i}}^{-\infty})
\end{equation*}
Here, $p_{m/2} ( \xi _{m/2} ) = 0$, and
$\theta _{i} = \arctan \frac{|\hat{F}_{n_{i}}(k_{n_{i}})|}{\epsilon _{n_{i}} } $. An
equivalent way to state the condition is
\begin{equation*}
\limsup \frac{- \log d( \cos ^{2} (\theta _{i} ), \Xi _{m} ) }{ \log N_{n_{i}} } = + \infty , \forall m\in \N ^{*} 
\end{equation*}
where $\Xi _{m} = \{ \xi ^{j} _{m} \}_{j}$ are the roots of the Legendre polynomial corresponding to
$m$.
\end{cond}

The $G_{\d }$ character of the condition is clear, since each $\theta _{i}$ depends
continuously on the cocycle in the $C^{\infty}$ topology. Density is also clear,
since by activating only higher modes one can make the
perturbation arbitrarily small in the $C^{\infty} $ topology. We make these arguments
precise in \S \ref{sec cond DUE}.

The lemma examines the effect of the change of coordinates on the cohomological
equation in the presence of obstructions. In the initial coordinates, we cannot obtain
an approximate solution with good estimates due to the resonance relation. The
problem is resolved by applying the change of coordinates by performing two steps
of the KAM scheme and then expressing the solution back in the initial coordinates.
If we follow the procedure of changing coordinates, solving the transformed equation
in the new coordinates, assuming condition \ref{condition DUE}, and undoing the
change of coordinates, we see that the Almost Reducibility and the corresponding
estimates allows us to conclude that we can eliminate the corresponding obstructions.
This is the content of the lemma which we now prove.

\begin{proof}[Proof of lemma \ref{iterative lemma}]
By lemma \ref{lem cohom obs}, in coordinates where $A_{1} = \{ e^{2i\pi (k_{1}\a +\epsilon _{1})} ,0 \} $ is diagonal,
the obstructions to the solution of the cohomological eq. over $(\a ,A_{1})$ are the modes
$e^{2i\pi k_{p} x}\p ^{m}_{j,p} (z,w)$ such that
\begin{equation} \label{def resonant modes}
k_{p}+(m-2p)k_{1}=0
\end{equation}
The operator $Ob_{1} = Ob^{k_{1}}_{N} $ is precisely the
projection on these modes, which by hypothesis are also the spectral
support of $\f  $. This implies that the equation
\begin{equation*}
\psi (x+\a , e^{-F_{1} (x)}A_{1}^{*}. \{ z,w \} )- \psi (x), \{ z,w \} ) =  \f  (x, \{ z,w \} )
\end{equation*}
cannot be solved with good estimates by direct application of  lemma \ref{lem cohom res 1 mode}.

We now prove the first point.
In the same coordinate system, the conjugation $B_{1} (\. )$ is of the form
\begin{equation*}
\begin{pmatrix}
e^{-2i\pi k_{1} \. } & 0 \\
0 & e^{2i\pi k_{1} \. }
\end{pmatrix}
\end{equation*}
And, by hypothesis,
\begin{equation*}
B_{1} (x+\a ). A_{1}.e^{\tilde{F}_{1} (x)}. B_{1}^{* } (x ) = \{ e^{2i\pi \epsilon _{1}} ,0 \}.
e^{\{ 0, \hat{F}_{1,z} (k_{1} \})} = A_{2}
\end{equation*}
By lemma \ref{transformation of harmonics}, such a conjugation
transforms the obstruction of the type
$e^{2i\pi k_{1} x}\p _{m}^{j,p}(z,w)$ to the harmonic $\p _{m}^{j,p}(z,w)$,
which proves the first point.

We now move on to proving the second point.
We remark that the vectors $\z ^{l} \w ^{l}$ are fixed by the reduction of resonant modes,
and so is the corresponding line in the matrix of $\p _{m}$, the harmonics $\p _{m}^{j,l}$ with $l = m/2$.
These are the only obstructions corresponding to $k=0$ in the 
coordinates $(x, \{ z,w \})$. Therefore, in the coordinates $(x, B_{1}(\. ). \{ z,w \})$ the rhs
of the cohomological equation is written in the form
\begin{eqnarray*}
\tilde{\f } (x,B_{1}(\. ). \{ z,w \} )  =
\widetilde {Ob}_{1} \f  (x,B_{1}(\. ). \{ z,w \} ) &=&
\sum _{j,p} \f  _{k_{p},j,p}^{m} \p _{m}^{j,p} (z,w) \\
&=&
\sum _{j,p} \tilde{\f } _{0,j,p}^{m} \p _{m}^{j,p} (z,w)
\end{eqnarray*}
where the summation is over the resonant modes as defined in eq. \ref{def resonant modes}. The
obstruction is now constant in $x$.


Let us we denote by $(x, \{ z_{2} , w_{2} \} ) $ the new coordinates in which $A_{2} $ is diagonal,
obtained by applying the conjugation $D_{1}$ in the notation of lemma \ref{conj in fin reg iter}.
In these coordinates, the linearized equation assumes the form
\begin{equation*}
\tilde{\psi} ' (x+\a , A_{2}^{*}. \{ z_{2} , w_{2} \} )- \tilde{\psi} '(x, \{ z_{2} , w_{2} \} ) = 
\sum _{j,p} \tilde{\f } ^{m,(1) }_{0,j,p} \p _{m}^{j,p} (z_{2},w_{2})
\end{equation*}
where a priori all harmonics in $\EE _{\p_{m}}$ are authorized in the summation. As we have seen, the obstructions to this type
of equation which are constant in $x$ are the harmonics of the form $\p _{m}^{j,m/2}$, i.e.
\begeq
Ob_{2}\left(
\sum _{j,p} \tilde{\f } ^{m,(1) }_{0,j,p} \p _{m}^{j,p} (z_{2},w_{2})
\right) =
\sum _{j} \tilde{\f } ^{m,(1) }_{0,j,m/2} \p _{m}^{j,m/2} (z_{2},w_{2})
\endeq
This proves the second point, since we apply isometries and projections with respect
to the $L^{2}$ norm.

We now prove the third point.
Applying $D_{1}^{*}$, in order to go back to coordinates where $A_{1}$ is diagonal,
will give in general a full matrix $\p _{m}^{j,p}$, $0 \leq j,p \leq m$. Then,
composition  with $(0 , B_{1}^{*} (\. ) )$ implies the third point, by lemma
\ref{transformation of harmonics}.

Let us now come to the fourth point. The obstructions $Ob _{2} (\widetilde{Ob} _{1} \f )$ are constant
in $x$, and remain so after conjugation by $(0 , B_{2} (\. ) )$, since $B_{2} (\. )$ commutes with $A_{2}$.
Application of $Ob_{3}$ gives a function satisfying
\begeq
\| Ob _{3} (\widetilde{Ob} _{2} \f ) \| _{L^{2}} \leq C N_{2}^{-s_{0}} \| \f  \| _{L^{2}}
\endeq
Conjugation by $(0,B^{*}_{1}.B^{*}_{2}(\. ))$ gives a function with (eventually) full support in harmonics,
$\p _{m}^{j,p}$, $0 \leq j,p \leq m$, and support in $[-m(N_{1}+N_{2}),m(N_{1}+N_{2})]$ in the frequency space,
which implies the fourth point, since $N_{1}+N_{2} \approx N_{2}$.
\end{proof}

\subsection{Proof of theorem \ref{thm DUE}}
We can now prove the existence of an approximate solution under condition \ref{condition DUE}.
Genericity of the condition is proved in the next section.
\begin{proof} [Proof of theorem \ref{thm DUE}]
It is just an iteration of the lemma \ref{iterative lemma} into the following algorithm,
where we call $\Phi = (\a , A.e^{F(\. )})$.

Fix $\f  \in \sm$ and $M$ big enough so that
\begeq
\f ' = \f  - \sum _{\substack{k \in \Z \\
0 \leq m \leq M \\
j,p
}} \f  _{k,j,p}^{m} \p _{m}^{j,p} (z,w)
\endeq
is small. Fix $i_{0}$ so that for $i \geq i_{0}$, $M N _{i} \ll N_{i+1}$.
Lemma \ref{iterative lemma} can now be applied iteratively for each $0 <m<M$.
For simplicity in notation, we will omit the prime in the notation, and denote the
truncated function by $\f = \f _{0}$ and implement the following algorithm.

At the first step, lemma \ref{lem cohom eq w pert} is applied to the function
\begeq
\f _{0}' = \f _{0} \circ (0,H_{i_{0}})
\endeq
and to the cocycle $\Phi _{i_{0}} = Conj _{H_{i_{0}}}\Phi $. Call $\psi _{0}'$ the solution,
transform it by $(0,H_{i_{0}})^{-1}$ in order to obtain
\begeq
\psi _{0} = \psi _{0}' \circ (0,H_{i_{0}})
\endeq
and call
\begeq
\f _{1} = \f _{0} - (\psi _{0} \circ \Phi - \psi _{0})
\endeq

In the next step, apply lemma \ref{lem cohom eq w pert} to the function $\f _{1}$ and
the cocycle $Conj _{H_{i_{0}+1}}\Phi $, and iterate.

The solution of the cohomological equation in step $i \geq i_{0}$ of the
K.A.M. scheme is carried out in the coordinates $(x, \{ z _{i} , w _{i} \} )$
where the cocycle is a perturbation of $(\a ,A_{i})$, and $A_{i}$ is diagonal.

At each step of the K.A.M. scheme, we solve the equation over the cocycle
$(\a ,A_{i})$. Lemma \ref{lem cohom eq w pert} provides a solution for the
non-resonant part of the function $\f$ expressed in the coordinates
$(x, \{ z _{i} , w _{i} \} )$.

Lemma \ref{lem cohom eq change coord} shows that
the solution for the non-resonnant part is settled, since the error term, denoted
by $E$ in the notation of the lemma, is of the order of $N_{n_{i}}^{-\infty }$, while
the corresponding change of coordinates is of the order of $N_{n_{i-1}}^{-\infty }$.

We thus need only solve for the obstructions, as they are defined in lemma
\ref{lem cohom obs}, applying iteratively lemma \ref{iterative lemma} in the following way.

The change of coordinates $(x,\{z,w \} ) \mapsto (x, \{ z _{i} , w _{i} \} ) $
is given $\prod _{1}^{i-1} (0,D_{i} )  (0,B_{i} (\. ) )$, where the $B_{i}(\. )$ are
produced by the K.A.M scheme, and the $D_{i}$, defined in $G / \NN _{i}$, diagonalizes
$A_{i}$ in the standard system of coordinates. Condition \ref{condition DUE}
implies that for any given $m =2l \in \N $, there exists a subsequence $\{ i_{j} \}$ along
which eq. \ref{Legendre lemma} holds for arbitrarily big $s_{0}$. This implies that the
error term becomes arbitrarily small in the standard coordinates, and the proof is
concluded.
\end{proof}

\section{Genericity of condition \ref{condition DUE}} \label{sec cond DUE}
By invoking lemma \ref{lem DUE Gd}, we can already conclude the proof of
theorem \ref{thm DUE}. We prefer, however, to show that condition \ref{condition DUE}
is itself generic, in order to make clearer the fact that DUE cocycles are obtained via
an Anosov-Katok construction, fibered over the rotation $\a$, in an otherwise quite
standard manner, which echoes in the genericity of the property that we establish.
Genericity is obtained roughly in the same way as genericity of condition of
eq. \eqref{condition UE}.

\begin{prop}
Condition \ref{condition DUE} is generic in $\WW$.
\end{prop}

\begin{proof}
Let $(\a , A)$ be a constant cocycle. It can be perturbed to the resonant constant cocycle $(\a , A_{1})$,
the resonant mode $k_{1}$ being arbitrarily big. Chose, now, $\xi _{m_{1}}^{j_{1}}$,
a root of a Legendre polymomial.
Then condition \ref{condition DUE} can be satisfied exactly for arbitrarly small
$\epsilon_{1}$ and $\hat{F}_{1}(k_{1})$. The condition is by definition open.

The cocycle thus constructed is reducible to $(\a , A_{2})=(\a ,\{ e^{2i\pi a_{2}},0\})$ where
\begeq
a_{2} = \sqrt{\epsilon_{1}^{2}+ |\hat{F}_{1}(k_{1})|^{2}}
\endeq
Thus, $A_{2}$ can be made resonant without changing $\theta _{1} =\arctan \frac{|\hat{F}_{1}(k_{1})|}{\epsilon_{1}}$.

This construction can be iterated along subsequences with different choices of $m_{i}$ in
order to satisfy condition \ref{condition DUE}, which is consequently generic by virtue of
theorem \ref{finite reduc thm}.
%
%
%
\end{proof}

In particular, the elimination of the modes $\p _{2} ^{1,1}$ (the line in the matrix of $\p _{2}$ corresponding
to the polynomial $\z \w$) demands that
\begin{equation*}
|\frac{|\hat{F}_{i_{l}^{(2)}} (k_{i_{l}^{(2)}} )|}{\epsilon _{n_{i_{l}^{(2)}}}} -1|
= O (N_{i_{l}^{(2)}}^{- \infty})
\end{equation*}
which implies a fortiori the condition of eq. \eqref{condition UE}.

\section{Proof of proposition \ref{thm cohom stab implies reduc}} \label{proof of thm 1.4}

In this section, we prove that cohomologically stable cocycles are
reducible, and consequently induce a Diophantine
rotation on their invariant tori.

Suppose that the cocycle $(\a , Ae^{F(\. )})$ is given in normal form
(for a fixed choice of parameters),
with an infinite number of steps $\{ n_{i} \}$ where resonances are reduced.
Let us first relate the size of the quantity $\epsilon _{i}$, i.e. the distance of the
cocycle at the $i$-th step of the KAM scheme from the exact resonance with the relevant
scale at that step. This amounts to calculating
\begin{equation} \label{eq def nbar}
\bar{\n } = \limsup \nu _{i} = \limsup  \frac{-\log |\epsilon _{i} |}{ \phantom{- }  \log N_{n_{i}}} \in [\n , \infty ]
\end{equation}
which tells us whether this distance is polynomial or super-polynomial with respect to the
order of truncation.

In the spirit of \cite{NKRigidity}, we can prove the following lemma.
\begin{lem} \label{lem finite nu}
If $\bar{\nu} < \infty$, then the cocycle is $\sm$ reducible.
\end{lem}
\begin{proof}
If $ \bar{\n } < \infty $, then the K.A.M. scheme with $\n ' =2 \bar{\n }$ has only a finite
number of resonant steps.

Such a scheme can be applied in the following way. We initiate the scheme with the fixed parameter $\n$,
and go deep enough so that the smallness conditions of the scheme with the new parameter
$\n '$ are satisfied and so that
\begeq
| \epsilon _{i} |> N_{i}^{-\nu '}
\endeq
for the remaining steps.
It is a matter of a direct calculation to verify that the new scheme does not produce
any reductions of resonant modes.

Alternatively, one can apply the conjugations constructed in lemma
\ref{conj in fin reg iter} and obtain a converging product of conjugations, since under
the assumption that $\bar{\nu} < \infty$ we have
\begeq
\theta _{i} = O(N_{i}^{-\infty})
\endeq
This concludes the proof.
\end{proof}

This procedure amounts to reducing the initial perturbation to such a size that
$\a $ can be placed in $DC (\gamma,  \t ')$, for a convenient $\t ' > \t$, and re-initiating the scheme for
the rotation satisfying this weaker Diophantine condition (see also
\cite{NKRigidity} for a similar argument). Anyhow, the cocycle in
this case is reducible and therefore not DUE. It is cohomologically
stable iff the induced rotation on a (and thus on each)
invariant torus $\T \times \Sp ^{1}$ is Diophantine.

In view of the previous lemma, in the proof of proposition \ref{thm cohom stab implies reduc}
we can assume that $ \bar{\n } = \infty $. This is weaker than condition
\ref{condition DUE}.
In this case, an infinite subsequence of steps of the scheme, still denoted by
$\{ n_{i} \}$, reduces essential resonances.
The argument in this case follows the proof of the
fact that Liouvillean rotations are not cohomologically stable, and this is relevant, since
$\bar{\nu} = \infty$ implies that the cocycle at the $i$-th step is
$O_{\sm }(N_{n_{i}}^{-\infty})$
close to a resonant cocycle, which makes its behaviour Liouvillean. 

\begin{proof}[Proof of proposition \ref{thm cohom stab implies reduc}]
By lemma \ref{lem finite nu}, we can assume that $ \bar{\n } = \infty $  and suppose that
such a cocycle is Cohomologically Stable.

We fix $m>0$ and a subsequence of $n_{i}$ such that
$ \nu _{i} \ra \infty $. At the step $i$, choose
$\f  ^{\prime} _{i} (x, \{ z_{i} , w_{i} \}) $ as in lemma \ref{lem cohom obs},
i.e. a function supported in the obstructions of the $i$-th step. Such functions form an
$(m+1)^{2}$-dimensional space.

By lemma \ref{transformation of harmonics}, as we have already observed in the proof of
lemma \ref{iterative lemma}, the dimension of the space of obstructions at the $i$-th step
that are transformed into obstructions at the $i+1$-th step is $m+1$. There exist,
therefore, $\f _{1}' \neq 0$ such that $\| Ob_{i+1}\f  ^{\prime} _{i} \| _{L^{2}} = 0$.
Consequently, there exists $\tilde{\psi} _{i}$ solving
\begin{equation*}
\tilde{\psi}  _{i} \circ (\a , A_{i+1}) - 
\tilde{\psi}  _{i} = \f ^{\prime} _{i}
\end{equation*}
satisfying, by lemma \ref{lem cohom res 1 mode},
\begin{equation}
\| \tilde{\psi}  _{i} \| _{H^{s}} \approx 
\r _{i} ^{-1 } N_{n_{i}}^{s+2}
\| \f ^{\prime} _{i}\|_{L^{2}}
\end{equation}
where
\begeq
\r _{i} = d (A_{i+1},\Id )
\endeq
as defined in eq. \ref{eq def rho} is $O(N_{i}^{-\infty})$ under the
assumption that $\bar{\nu} = \infty$.

Therefore, $\f  _{i}^{\prime}$ is $O_{C^{\infty}}(N_{i+1}^{-\infty})$
close to a coboundary over $(\a , A_{i+1} e^{F_{i+1}})$, since
\begin{equation*}
\tilde{\f } _{i} = \f ^{\prime} _{i} +
\tilde{\psi} _{i} \circ (\a , A_{i+1}
e^{F_{i+1}})
- \tilde{\psi}  _{i} \circ (\a , A_{i+1})
\end{equation*}
is a coboundary, and satisfies
\begin{equation*}
\| \tilde{\f } _{i} - \f ^{\prime} _{i} \| _{H^{s}} \approx
\|F_{i+1} \|_{s}\r _{i}^{-1} N_{n_{i}}^{s+4}
\| \f ^{\prime} _{i}\|_{L^{2}}
\end{equation*}
We remind that, since $ \bar{\n } = \infty $,
eq. \ref{rareness of red of res} implies that $n_{i+1}-n_{i} \ra \infty$
(where one could replace $\n $ by $\n _{i} \ra \infty$), which in turn implies
the estimate.

We now impose the condition
\begeq
\|\f  ^{\prime} _{i} \|_{L^{2}} = \r _{i} = O(N_{i}^{-\infty})
\endeq
We also call $\f _{i}$ the function $\tilde{\f } _{i}$ expressed
in the standard coordinates, and call
\begin{equation*}
\f  (x , \{ z ,w \}) = \sum _{i =0 }^{\infty} \f _{i} (x , \{ z ,w \})
\in C_{0} ^{\infty} (\T \times G)
\end{equation*}
The function $\f  $ is a coboundary, since it is a limit of coboundaries over
$(\a , Ae^{F(\. )})$, and this latter space is
closed by hypothesis. Let $\psi \in C ^{\infty} $ be the solution of the cohomological equation:
\begin{equation*}
(\a , Ae^{F(\. )}).\psi -  \psi = \f 
\end{equation*}
If we apply the change of coordinates $H_{n}(\. )$ to this equation, we obtain
\begin{equation*}
(\a , A_{i}e^{F_{i}(\. )}).\psi \circ H_{i} -  \psi \circ H_{i}= \f  \circ H_{i}
\end{equation*}
which implies directly that
\begin{equation*}
(\a , A_{i}).\psi \circ H_{i} -  \psi \circ H_{i}= \f  \circ H_{n} + O_{\sm}(N_{i}^{-\infty})
\end{equation*}

We now project this equation on the shell of frequencies
$[\frac{1}{2}N_{n_{i}} , \frac{3}{2} N_{n_{i}} ] $. Then
\begeq
T_{\frac{1}{2}N_{n_{i}}}^{\frac{3}{2}N_{n_{i}}} \psi =
R_{\frac{1}{2}N_{n_{i}}}T_{\frac{3}{2}N_{n_{i}}} \psi
\endeq
must carry an $L^{2}$ norm of the order of constants, since
\begeq
T_{\frac{1}{2}N_{n_{i}}}^{\frac{3}{2}N_{n_{i}}}(\psi \circ H_{i} )\approx \psi '_{i}
\endeq
as is satisfies
\begin{equation*}
(\a , A_{i}).T_{\frac{1}{2}N_{n_{i}}}^{\frac{3}{2}N_{n_{i}}}(\psi \circ H_{i}) -
T_{\frac{1}{2}N_{n_{i}}}^{\frac{3}{2}N_{n_{i}}}(\psi \circ H_{i})= \tilde{\f} _{i} 
+ O_{\sm}(N_{i}^{-\infty})
\end{equation*}

This is incompatible with the hypothesis that $\psi \in C^{\infty}$.
\end{proof}

\section{Necessity of condition \ref{condition DUE}}

In this section we let condition \ref{condition DUE} degenerate in a
controlled way in order to show that a condition of this kind is indeed necessary in
order to obtain DUE. We prove the following proposition, which is to be compared with
section \ref{loc dens of fin dif}.
\begin{prop} \label{optim of cdtion}
Let $m \in \N$ be fixed, and $\Xi _{m} = \{ \xi ^{j} _{m} \}_{j}$ the roots of the corresponding Legendre polynomial. Let also
$(\a ,Ae^{F(\. )})$ be in normal form. Suppose, moreover, that
\begin{equation*}
\limsup \frac{-\log d( \cos ^{2}(\z _{i} )  , \Xi _{m} ) }{\log N_{n_{i}}} = s_{0} < \infty
\end{equation*}
for some $ \l < s_{0 } < \infty $ and
\begin{equation*}
\limsup \frac{-\log d ( \cos ^{2} (\z _{i} )  , \Xi _{m'} ) }{\log N_{n_{i}}} = \infty , m' \neq m
\end{equation*}
Then, $C^{\infty}$ coboundaries are dense in $C^{s_{0} - \l }$ but not in $C^{s_{0} + \l }$.
\end{prop}
\begin{proof}
The proof that $C^{\infty}$ coboundaries are dense in
$C^{s_{0} - \l }$ is just a finite regularity vesion
of the corresponding part of the proof of theorem \ref{thm DUE}, in \S
\ref{sec constr approx sol}.
The proof that $C^{\infty}$ coboundaries are not dense in $C^{s_{0} + \l }$ goes like the proof
of proposition \ref{thm cohom stab implies reduc} in \S \ref{proof of thm 1.4}.
\end{proof}

\section{Conclusions and remarks}

In this paper we made systematic use of the fact that every cocycle in the KAM regime
$\WW$ is Almost Reducibile, and moreover with good control over the size of the conjugation
with respect to the size of the perturbation at each step. This allowed us to sharpen the
conclusions of \cite{El2002a}, where this property was exploited in order to establish
ergodicity with respect to the Haar measure.

The technique that we applied can be described as a fibered Anosov-Katok construction, where
the rotation in the basis, $\T ^{d}$, remains fixed, and the approximation-by-conjugation
argument is carried out in the fibers, in an otherwise classical manner. The role of
periodic rotations is played by resonant cocycles, and Almost Reducibility guarantees that
the conclusions of the construction (i.e. genericity of DUE) hold in an open set.

The same point of view was used by the author in \cite{NKContSpec}, where a spectral
dichotomy was obtained for cocycles in $\WW$, and abstracted from in \cite{NKCohomRigAK},
where the intuitively obvious statement that Anosov-Katok constructions produce
Liouvillean, and therefore not Cohomologically Rigid, dynamics, is given a precise content.

\bibliography{aomsample}
\bibliographystyle{aomalpha}

\end{document}